\documentclass[a4paper,reqno,twoside]{amsart}

\usepackage{amssymb}
\usepackage{latexsym}
\usepackage{amsmath}
\usepackage{mathrsfs}
\usepackage{pifont}
\usepackage{enumerate}

\theoremstyle{plain}
\newtheorem{propn}{Proposition}[section]
\newtheorem{thm}[propn]{Theorem}
\newtheorem{lemma}[propn]{Lemma}
\newtheorem{cor}[propn]{Corollary}

\theoremstyle{definition}
\newtheorem*{defn}{Definition}

\theoremstyle{remark}
\newtheorem*{rem}{Remark}
\newtheorem*{rems}{Remarks}
\newtheorem*{notn}{Notation}

\theoremstyle{definition}

\newcommand{\ve}{\varepsilon}

\newcommand{\ida}{1_{\Bil}}

\newcommand{\Al}{\mathsf{A}}
\newcommand{\Bil}{\mathsf{B}}
\newcommand{\blg}{\mathsf{B}}
\newcommand{\alg}{\mathsf{A}}
\newcommand{\Cil}{\mathsf{C}}
\newcommand{\El}{\mathsf{E}}
\newcommand{\Eil}{\mathsf{E}}
\newcommand{\elg}{\mathsf{E}}
\newcommand{\Hil}{\mathsf{H}}
\newcommand{\hil}{\mathsf{h}}
\newcommand{\Kil}{\mathsf{K}}
\newcommand{\kil}{\mathsf{k}}
\newcommand{\Uil}{\mathsf{U}}
\newcommand{\Vil}{\mathsf{V}}
\newcommand{\Wil}{\mathsf{W}}
\newcommand{\Xil}{\mathsf{X}}
\newcommand{\Yil}{\mathsf{Y}}

\newcommand{\init}{\mathfrak{h}}
\newcommand{\noise}{\mathsf{k}}

\newcommand{\Ao}{\mathsf{A}_0}

\newcommand{\bc} {\mathbb C}
\newcommand{\bn}{\mathbb N}
\newcommand{\br}{\mathbb R}

\newcommand{\Alg}{\mathcal{A}}

\newcommand{\Blg}{\mathcal{B}}
\newcommand{\Vlg}{V}

\newcommand{\Exps}{\mathcal{E}}
\newcommand{\FFock}{\mathcal{F}}
\newcommand{\FFocktot}{\FFock_{[0,t[}}

\newcommand{\Op}{\mathcal{O}}
\newcommand{\Opdagger}{\mathcal{O}^\ddagger}

\newcommand{\Skorohod}{\mathcal{S}}

\newcommand{\ExpsD}{\mathcal{E}_D}

\newcommand{\Prob}{\mathbb{P}}

\newcommand{\Proc}{\mathbb{P}}
\newcommand{\cProc}{\mathbb{P}_{\fullstar}}

\newcommand{\cProcdagger}{\mathbb{P}^\ddagger_{\fullstar}}
\newcommand{\cProccb}{\mathbb{P}_{\fullstar\,\mathrm{cb}}}

\newcommand{\cR}{R_{\fullstar}}

\newcommand{\Pto}{\!\to\!}

\newcommand{\Procdagger}{\mathbb{P}^\ddagger}
\newcommand{\Proccb}{\mathbb{P}_{\mathrm{cb}}}

\newcommand{\Step}{\mathbb{S}}

\newcommand{\fullstar}{\star}
\newcommand{\Com}{\Delta}
\newcommand{\Cou}{\counit}

\newcommand{\counit}{\epsilon}
\newcommand{\CBcoproduct}{CB^{\Com}}
\newcommand{\Projk}{\Delta\!^{QS}}
\newcommand{\dQS}{\Projk}

\newcommand{\otM}{\otimes_{\Mat}}
\newcommand{\fullcomp}{\bullet}

\newcommand{\Real}{\mathbb{R}}
\newcommand{\Rplus}{\Real_+}
\newcommand{\Comp}{\mathbb{C}}
\newcommand{\Nat}{\mathbb{N}}

\newcommand{\ZZ}{\mathbb{Z}}

\newcommand{\Mat}{\mathrm{M}}

\newcommand{\cb}{{\text{\tu{cb}}}}

\newcommand{\wh}{\widehat}
\newcommand{\wt}{\widetilde}
\newcommand{\ol}{\overline}
\newcommand{\ot}{\otimes}

\newcommand{\olot}{\overline{\otimes}}

\newcommand{\op}{\oplus}
\newcommand{\la}{\langle}
\newcommand{\ra}{\rangle}

\newcommand{\chat}{\wh{c}}

\newcommand{\Dhat}{\wh{D}}
\newcommand{\fhat}{\wh{f}}
\newcommand{\ghat}{\wh{g}}

\newcommand{\kilhat}{\wh{\kil}}

\newcommand{\ktilde}{\wt{k}}

\newcommand{\tu}{\textup}

\DeclareMathOperator{\Dom}{Dom}
\DeclareMathOperator{\Ran}{Ran}
\DeclareMathOperator{\Lin}{Lin}
\DeclareMathOperator{\Ker}{Ker}
\DeclareMathOperator{\opposite}{op}
\DeclareMathOperator{\id}{id}
\DeclareMathOperator{\re}{Re}

\newenvironment{alist}
{

\begin{enumerate}}
{\end{enumerate}}

\newenvironment{rlist}
{

\begin{enumerate}}
{\end{enumerate}}

\numberwithin{equation}{section}
\begin{document}

\title
[Quantum stochastic convolution cocycles]
{Quantum stochastic convolution cocycles II}

\author[Lindsay]{J.\ Martin Lindsay}
\author[Skalski]{Adam G.\ Skalski}
\footnote{\emph{Permanent address of AGS}.
Department of Mathematics, University of \L\'{o}d\'{z}, ul.
Banacha 22, 90-238 \L\'{o}d\'{z}, Poland.}
\address{
Department of Mathematics and Statistics,
Lancaster University, Lancaster LA1 4YF, United Kingdom}
\email{j.m.lindsay@lancaster.ac.uk}
\email{a.skalski@lancaster.ac.uk}


\begin{abstract}
Sch\"urmann's theory of quantum L\'evy processes, and more generally the
theory of quantum stochastic convolution cocycles, is extended to the
topological context of compact quantum groups and operator space
coalgebras.
Quantum stochastic convolution cocycles on a $C^*$-hyperbialgebra, which
are Markov-regular, completely positive and contractive, are shown to
satisfy coalgebraic quantum stochastic differential equations with
completely bounded coefficients, and the structure of their stochastic
generators is obtained. Automatic complete boundedness of a class of
derivations is established, leading to a characterisation of the stochastic
generators of *-homomorphic convolution cocycles on a $C^*$-bialgebra.
Two tentative definitions of quantum L\'evy process on a compact quantum
group are given and, with respect to both of these, it is shown that an
equivalent process on Fock space may be reconstructed from the generator
of the quantum L\'evy process.
In the examples presented, connection to the algebraic theory is
emphasised by a focus on full compact quantum groups.
\end{abstract}

\maketitle

\section*{Introduction}

In this paper we investigate quantum stochastic evolutions with
independent identically distributed increments on compact quantum groups,
in other words \emph{quantum L\'evy processes}. The natural setting for
this analysis is the somewhat wider one of quantum stochastic convolution
cocycles. For a compact quantum group $\blg$, a quantum stochastic
convolution cocycle on $\blg$ is a family of linear maps $(l_t)_{t\geq 0}$
from $\blg$ to operators on the symmetric Fock space $\FFock$, over a
Hilbert space of the form $L^2(\Rplus; \noise)$, satisfying
\[
l_{s+t} = l_s \fullstar (\sigma_s \circ l_t),
\;\;\; s,t \geq 0
\]
and some regularity and natural adaptedness conditions. Here
$\big(\sigma_s\big)_{s\geq 0}$ is the semigroup of time-shifts on
$B(\FFock)$ and the convolution is induced by the quantum group
structure; the initial condition is specified by the counit: $l_0 =
\counit ( \cdot ) I_\FFock$.
Thus the increment of the process over the interval $[0,s+t]$ coincides
with the increment over $[0,s]$ convolved with the (shifted) increment
over $[0,t]$. We show that such families may be obtained as solutions of
quantum stochastic differential equations with completely bounded
coefficients, we analyse their positivity and multiplicativity properties,
and we establish natural conditions under which all sufficiently regular
cocycles arise in this way.
Motivated by these results (and the purely algebraic theory), we propose
two abstract definitions of quantum L\'evy process on a compact quantum
group and show that any process which has bounded `generator' has an
equivalent Fock space realisation. Precise definitions are given below.

Stochastic cocycles on operator algebras were introduced by Accardi
(under the name quantum Markovian cocycles) for
Feynman-Kac type perturbation of quantum dynamical semigroups
(\cite{acc}). Earlier work on a cocycle approach to classical Markov
processes and their It\^o integral representation may be found in
\cite{Pinsky}. Quantum stochastic differential equations (\cite{HuP})
were quickly seen to provide examples of stochastic cocycles and in fact
to characterise large classes of them in the Fock space context (see
\cite{Lgreifswald} and references therein).

The theory of quantum L\'evy processes, developed by Sch\"urmann and
others, generalises the classical theory of L\'evy processes on groups
(\cite{Heyer}), and Skorohod's theory of stochastic semigroups
(\cite{sko}), to the context of quantum groups or, more generally,
*-bialgebras (see \cite{Schurmann}, \cite{FranzSchott}, \cite{glo} and
references therein). A quantum L\'evy process on a quantum group $\Blg$ is
a time-indexed family of unital *-homomorphisms from $\Blg$ to some
noncommutative probability space, with identically distributed and
(tensor-)independent increments, satisfying the convolution increment
relation given by the coproduct of $\Blg$, and with initial condition
given by the counit of $\Blg$.
Sch\"urmann showed that each quantum L\'evy process may be equivalently
realised in a symmetric Fock space as a solution of a quantum stochastic
differential equation.
This led us to introduce and investigate, in this algebraic context, quantum
stochastic convolution cocycles (\cite{LSqscc1}). These are
linear (but not necessarily unital or *-homomorphic) maps from a coalgebra
to a space of Fock space operator processes, satisfying the convolution
increment relation and counital initial condition.

In the last twenty years there has been a growing interest in the theory
of \emph{topological} quantum groups. Starting from the fundamental
paper of Woronowicz (\cite{wor1}), where the concept of compact quantum
groups was first introduced (under the name of compact matrix
pseudogroups), it has led to a rich and well-developed theory, with a
satisfactory notion of \emph{locally compact} quantum group eventually
emerging in the work of Kustermans and Vaes (\cite{kuv}).
The main object becomes a $C^*$-algebra, equipped with a coproduct and
counit satisfying a corresponding form of coassociativity and counit
relations.

In this paper we go beyond the purely algebraic context treated
in~\cite{LSqscc1} and initiate the
study of quantum L\'evy processes on a compact quantum group, or more
generally on a $C^*$-bialgebra. Heeding P.-A.\,Meyer's dictum once more,
we again set our work in the wider context of quantum stochastic
convolution cocycles on a coalgebra.
The coalgebras here though are operator-space-theoretic rather than being
purely algebraic. Nevertheless the stochastic cocycles in question may be
obtained by solving coalgebraic quantum stochastic differential equations.
In turn, every sufficiently regular completely positive and contractive
quantum stochastic convolution cocycle on a $C^*$-hyperbialgebra is shown
to satisfy a quantum stochastic differential equation of the above type.
These results are obtained by, on the one hand applying
techniques of operator space theory (\cite{EfR}, \cite{Pisier}), and on
the other hand using known facts about standard quantum stochastic
cocycles (see~\cite{LWjfa}, \cite{Lgreifswald} and references therein).
Here it is natural to work with processes on abstract operator spaces and
$C^*$-bialgebras. For this we use theory developed in~\cite{Sanhan} and
summarised in the first section. When the spaces are concrete this reduces
to the existing theory. A key tool of our analysis is a convolution
operation which we call the $R$-\emph{map}. This transforms coalgebraic
objects such as convolution cocycles and coalgebraic quantum stochastic
differential equations to standard objects of quantum stochastic analysis
(\cite{Lgreifswald}), setting up a traffic of properties and relationships
which we systematically exploit.
The $R$-map gives rise to a noncommutative avatar of the transformation
between convolution semigroups of measures and Markov semigroups (of
operators), familiar from classical probability theory.

The structure of the stochastic generators of Markov-regular,
*-homomorphic convolution cocycles on a $C^*$-bialgebra may be
characterised in terms of $\counit$-\emph{structure maps}, where $\counit$
is the counit,  or \emph{topological} Sch\"urmann triples (cf.\ their
purely algebraic counterparts). The complete boundedness of such
generators, indeed their implementability, follows from their
algebraic properties alone. We prove this by first extending well-known
results of Sakai, Ringrose and Christensen, on automatic continuity and
innerness properties of derivations, to the case of
$(\pi',\pi)$-derivations.  The fact that every $\counit$-structure map
defined on the whole $C^*$-bialgebra must  be implemented may be viewed as
a noncommutative counterpart to the fact that every classical L\'evy
process on a topological group which has a bounded generator must be a
compound Poisson process. In this connection we note the definition of
quantum Poisson process on a *-bialgebra proposed in~\cite{UweGreifswald}.

The axiomatisation of quantum L\'evy processes on a $C^*$-bialgebra raises
several problems connected with the fact that the product on a
$C^*$-algebra $\alg$ usually fails to extend to a continuous map from the
spatial tensor product $\alg \ot \alg$ to the algebra. We offer two
different ways of overcoming this obstacle, for both of which a
topological version of Sch\"urmann's reconstruction theorem remains valid.

Our choice of examples is designed to expose the variety of connections
of this work with the classical and quantum probabilistic literature. The
analysis of quantum stochastic convolution cocycles in the topological
context requires different methods and techniques to that of the purely
algebraic and poses new nontrivial problems.
However, according to our philosophy (explicitly described in the
expository paper \cite{LSbedlewo}),  purely algebraic and topological
convolution cocycles may nevertheless usefully be viewed from a common
vantage point.
This perspective is particularly well illustrated in the last class of
examples discussed here, namely that of *-homomorphic quantum stochastic
convolution cocycles on a full compact quantum group. We would also like
to point out that conversely, due to the Fundamental Theorem on
Coalgebras, one can
view the purely algebraic situation as a \emph{finite dimensional version}
of the topological theory.
An example of reasoning along these lines may be found in the final
section of \cite{LSqsde}.

The plan of the paper is as follows. In the first section we review the
basic facts needed from operator space theory and quantum stochastic
analysis. We work with processes on abstract operator spaces. The
transition from concrete to abstract exploits a number of natural
identifications and inclusions, the key ones
being~\eqref{bilinear}, ~\eqref{sp incl cb} and~\eqref{orawan}.
In Section~\ref{sec: B} the notion of operator space coalgebra
is introduced and basic properties of the $R$-map are established,
facilitating a correspondence between mapping composition structures and
convolution-type structures.
Section~\ref{sec: D} contains proofs of the existence, uniqueness and
regularity of solutions of coalgebraic quantum stochastic differential
equations with completely bounded coefficients. There also the ground is
prepared for a traffic between standard quantum stochastic cocycles and
quantum stochastic convolution cocycles.
The latter are defined in Section~\ref{sec: E} where the solutions of
coalgebraic quantum stochastic differential equations are shown to lie in
this class. The section concludes with a
brief discussion of \emph{opposite} convolution cocycles.
In Section~\ref{sec: F} the converse result is established for
Markov-regular, completely positive, contractive quantum stochastic
convolution cocycles on a $C^*$-hyperbialgebra: they are characterised as
solutions of coalgebraic quantum stochastic differential equation with
completely bounded coefficient of a particular form.
Section~\ref{sec: G} deals with *-homomorphic convolution cocycles on a
$C^*$-bialgebra. As in the purely algebraic case, their stochastic
generators are characterised by structure relations involving the counit;
in the topological case these amount to the generator being an
$\counit$-structure map where $\counit$ is the counit of the bialgebra.
In Section~\ref{Levy} two candidates for the axiomatisation of quantum
L\'evy processes on a $C^*$-bialgebra are proposed; firstly, in a weak
sense of distributions, and secondly, as processes whose values are
operators from a product system, in the sense of Arveson.
Basic consequences of the proposed axioms are discussed, and
reconstruction theorems established.
Section~\ref{section: H} is devoted to examples, first the commutative
case of classical compact groups, then the cocommutative case of the
universal $C^*$-algebra of a discrete group, and finally the case of full
compact quantum groups. In the latter case a link is established with the
purely algebraic quantum stochastic convolution cocycles investigated in
\cite{LSqscc1}.
In an appendix some results on
derivations are established; these are applied
to yield the automatic implementedness of $\counit$-structure maps
used in Section~\ref{sec: G}.

Some of the results proved here have been announced in \cite{LSbedlewo}.

\medskip
\noindent
\emph{Note added in proof.}
It is now clear that our results extend to the context of locally compact
quantum groups in the sense of Kustermans and Vaes (\cite{LSqscc3}).

\subsection{Notation}
All vector spaces arising in this paper are complex; inner products (and
all sesquilinear maps) are linear in their \emph{second} argument.
For a dense subspace $E$ of a Hilbert space $\hil$,
$\Op (E)$ denotes the space of operators $\hil\to\hil$ with domain
$E$ and
$\Opdagger (E):= \{ T\in\Op (E): \Dom T^*\supset E \}$.
Thus $\Opdagger (E)$ has the natural
conjugation $T\mapsto T^{\dagger}:=T^*|_{E}$.
We view
$B(\hil)$ as a subspace of $\Opdagger (E)$ (via restriction/continuous
linear
extension).
For vectors $\zeta\in E$ and $\zeta'\in \hil$,
$\omega_{\zeta',\zeta}$ denotes the linear functional on
$\Op (E)$ given by
$T\mapsto \la\zeta', T\zeta\ra$.
We use the Dirac-inspired notations
\[
|E\ra := \{|\zeta\ra : \zeta \in E\}
\text{ and }
\la E| := \{\la\zeta | : \zeta \in E\}
\]
where $|\zeta\ra\in |\hil\ra := B(\Comp;\hil)$ and
$\la \zeta |\in \la\hil | := B(\hil;\Comp)$ are defined by
$\lambda \mapsto \lambda\zeta$ and $\eta\mapsto\la\zeta, \eta\ra$
respectively.
A class of ampliations frequently met here is denoted as follows:
\begin{equation}\label{ampliation}
\iota_{\hil}:\Vil\to\Vil\ot B(\hil), \quad x\mapsto x\ot I_{\hil},
\end{equation}
where this time the operator space $\Vil$ is determined by context (and
$\ot$ denotes spatial tensor product).

For a vector-valued function $f$ on $\Rplus$ and subinterval $I$ of
$\Rplus$ $f_I$ denotes the function on $\Rplus$ which agrees with $f$ on
$I$ and vanishes outside $I$. Similarly, for a vector $\xi$, $\xi_I$ is
defined by viewing $\xi$ as a constant function. This extends the standard
indicator function notation.
The symmetric measure space over the Lebesgue measure space $\Rplus$
(\cite{Guichardet})
is denoted $\Gamma$, with integration denoted
$\int_{\Gamma} \cdots d\sigma$, thus
$\Gamma = \{\sigma\subset\Rplus : \# \sigma < \infty\} =
\bigcup_{n\geq 0} \Gamma^{(n)}$ where
$\Gamma^{(n)} = \{\sigma\subset\Rplus : \# \sigma =n \}$
and $\emptyset$ is an atom having unit measure.
If $\Rplus$ is replaced by a subinterval $I$ then we write $\Gamma_I$ and
$\Gamma^{(n)}_I$,
thus the measure of $\Gamma^{(n)}_I$ is $|I|^n/n!$ where $|I|$ denotes the
length of $I$.

For a linear map $\psi : U\to V$ the corresponding linear map
between conjugate vector spaces
\begin{equation} \label{conjugate map}
U^\dagger \to V^\dagger,
\quad x^\dagger \mapsto \psi (x)^\dagger
\end{equation}
is denoted $\psi^\dagger$; $L(U^\dagger; V^\dagger)$ is thereby
the natural conjugate space of $L(U;V)$.
The collection of sesqilinear maps $\phi :U\times V\to W$
is denoted $SL(U,V;W)$;
when $W$ is a space of maps we denote values of $\phi$ by $\phi^{u,v}$
($u\in U, v\in V$). The collection of bilinear maps $U\times V\to W$ is
denoted $L(U,V;W)$.
If $\Alg$ is an involutive algebra and $E$ is a dense subspace of a
Hilbert space $\hil$ then \emph{weak multiplicativity} for a map $\phi
:\Alg \to \Opdagger (E)$, is the property
\begin{equation} \label{weak mult}
\phi(a^*b) = \phi^\dagger (a)^*\phi(b) \quad (a,b\in\Alg),
\end{equation}
where $\phi^\dagger : a\mapsto \phi(a^*)^*|_E$.

\begin{rem}
If $\phi : \Al \to \Opdagger (E)$ is a linear map, defined on a
$C^*$-algebra, which is \emph{real} (that is
$\phi = \phi^\dagger$) and weakly multiplicative then $\phi$ is
necessarily bounded-operator-valued and thus may be viewed as a *-homomorphism
$\Al\to B(\hil)$.
\end{rem}

\section{Operator space and quantum stochastic preliminaries} \label{sec: A}

In this section we collect some relevant facts from operator space theory,
recall the matrix-space construction
and describe the basic properties of tensor-extended compositions.
We also recall relevant results from quantum stochastic (QS) analysis.

\subsection{Operator spaces (\cite{EfR}, \cite{Pisier})}
For operator spaces $\Vil$ and $\Wil$
the Banach space of completely bounded maps from $\Vil$ to $\Wil$ is endowed with operator space
structure via the linear identifications
\[
M_n \big( CB(\Vil; \Wil) \big) =
CB\big(\Vil; \Mat_n (\Wil) \big) \quad (n \in \Nat) ,
\]
where $\Mat_n (\Wil)$ denotes the linear space
$M_n (\Wil)$ with its natural operator space structure.
When viewed as a $C^*$-algebra or operator space, $M_n(\Comp)$ is denoted
$\Mat_n$.

The operator space spatial/minimal tensor product of $\Vil$ and $\Wil$ is
here denoted simply $\Vil \ot \Wil$.
For example $\Mat_n(\Wil)$ may be identified with the spatial tensor
product $\Wil\ot\Mat_n$.
When $\Vil$ and $\Wil$ are realised in $B(\Hil)$ and $B(\Kil)$
respectively, $\Vil \ot \Wil$ is realised concretely in
$B(\Hil \ot \Kil) = B(\Hil) \ol{\ot} B(\Kil)$ as the norm closure of the
algebraic tensor product $\Vil \odot \Wil$.
In fact $\Vil \ot \Wil$ does not depend on concrete realisation of $\Vil$
and $\Wil$; an abstract model arises from the natural linear embedding
\begin{equation} \label{abstract model}
\Vil \odot \Wil \hookrightarrow CB(\Wil^* ;\Vil)
\end{equation}
(where $\Wil^*$ is defined below). For any completely bounded maps $\phi : \Vil \rightarrow \Vil'$ and $\psi
: \Wil \rightarrow \Wil'$ into further operator spaces, the linear map
$\phi \odot \psi$ extends uniquely to a completely bounded map $\Vil \ot
\Wil \rightarrow \Vil' \ot \Wil'$; the extension is denoted $\phi \ot
\psi$ and satisfies
$\| \phi \ot \psi\|_{\cb} = \| \phi \|_{\cb} \| \psi \|_{\cb}$.
Each bounded operator $\phi : \Vil \rightarrow \Mat_n$ is automatically
completely bounded and satisfies
$\| \phi \|_{\cb} = \| \phi^{(n)} \|$ in the notation
$\phi^{(n)}: [x_{ij}]\mapsto [\phi(x_{ij})]$, in other words
$\phi^{(n)}=\phi \ot \id_{\Mat_n}$. In particular, the operator space
$CB(\Vil;\Comp)$
coincides with the Banach space dual $B(\Vil; \Comp)$ and has the same norm; it is therefore
denoted $\Vil^*$.
Note the natural completely isometric isomorphisms
\begin{equation} \label{bilinear}
CB(\Uil,\Vil;\Wil) = CB\big(\Uil; CB(\Vil;\Wil)\big)
\end{equation}
for operator spaces $\Uil$, $\Vil$ and $\Wil$.
We shall also exploit the natural completely isometric inclusions
\begin{equation} \label{sp incl cb}
\Vil\ot B(\Hil;\Hil') \hookrightarrow CB\big(\la\Hil'|, |\Hil\ra;\Vil\big)
\end{equation}
for operator space $\Vil$ and Hilbert spaces $\Hil$ and $\Hil'$.
(See below for the tensor product which delivers \emph{isomorphism} here.)

The following short-hand notation for tensor-extended composition
is useful. Let $\Uil, \Vil , \Wil$ and $\Xil$ be operator spaces, and let
$V$ be a vector space. If $\phi \in L(V ; \Uil \ot \Vil \ot \Wil)$ and
$\psi \in CB(\Vil ; \Xil)$ then we compose in the obvious way:
\begin{equation} \label{fullcomp}
\psi \fullcomp \phi :=
(\id_{\Uil} \ot \psi \ot \id_{\Wil} ) \circ \phi
\in L(V;\Uil \ot \Xil \ot \Wil).
\end{equation}
Ambiguity is avoided provided that the context dictates which
tensor component the second-to-be-applied map $\psi$ should act on.
This also applies to the case where
$\phi\in
SL\big(\Hil',\Hil;L(V;\Vil)\big)
$
as follows:
$
\psi\fullcomp \phi \in
SL\big(\Hil',\Hil;L(V;\Xil)\big)
$
is given by
\begin{equation} \label{abstract full dot}
\big(\psi\fullcomp \phi\big)^{\xi',\xi} =
\psi\circ\phi^{\xi',\xi}.
\end{equation}
The natural inclusion
$L\big(V; CB\big( \la\Hil'|,|\Hil\ra;\Vil\big)\big)
\subset
SL\big(\Hil',\Hil;L(V;\Vil)\big)$
is relevant here.

\subsection{Matrix spaces (\cite{LWfeller})}
For an operator space $\Yil$ in $B(\Hil;\Hil')$ and Hilbert spaces $\hil$
and  $\hil'$ define
\begin{equation} \label{matrix space defined}
\Yil\otM B(\hil;\hil') :=
\{ T\in B(\Hil\ot\hil; \Hil'\ot\hil') = B(\Hil;\Hil') \olot B(\hil;\hil'):
\Omega^{\Yil}_{\zeta',\zeta}(T)\in\Yil \}
\end{equation}
where $\Omega^{\Yil}_{\zeta',\zeta}$ denotes the slice map
$\id_{\Yil} \olot\omega_{\zeta',\zeta}$.
For us the relevant cases are
$\Yil\otM B(\hil)$ and
$\Yil\otM |\hil\ra$, referred to respectively as the
$\hil$-matrix space over $\Yil$ and the
$\hil$-column space over $\Yil$.
Matrix spaces are operator spaces which lie between the spatial tensor
product $\Yil\ot B(\hil;\hil')$ and the ultraweak tensor product
$\ol{\Yil}\olot B(\hil;\hil')$, coinciding with the latter when $\Yil$ is
ultraweakly closed
($\ol{\Yil}$ here denotes the ultraweak closure of
$\Yil$).
They arise naturally in quantum stochastic
analysis where a topological state space is to be coupled with the
measure-theoretic noise --- if $\Yil$ is a $C^*$-algebra then typically
the inclusion $\Yil\ot B(\hil)\subset \Yil\otM B(\hil)$ is proper and
$\Yil\otM B(\hil)$ is \emph{not} a $C^*$-algebra.
Completely bounded maps between concrete operator spaces lift to
completely bounded maps between corresponding matrix spaces: if $\Yil'$ is
another concrete operator space,
for $\phi\in CB(\Yil;\Yil')$ there is a unique map
$\Phi : \Yil\otM B(\hil;\hil')\to\Yil'\otM B(\hil;\hil')$ satisfying
\[
\Omega^{\Yil'}_{\zeta',\zeta}\circ \Phi = \phi\circ\Omega^{\Yil}_{\zeta',\zeta}
\quad
(\zeta\in\hil, \zeta\in\hil');
\]
it is denoted $\phi \otM \id_{B(\hil;\hil')}$.
Using these matrix liftings, tensor-extended
compositions work in the same way for matrix spaces as for spatial tensor products.
There are natural completely isometric isomorphisms
\begin{equation} \label{orawan}
\Yil\otM B(\hil;\hil') = CB(\la\hil' |, |\hil\ra;\Yil)
\end{equation}
(cf.\ \eqref{sp incl cb}) under which
$\phi \otM \id_{B(\hil;\hil')}$
corresponds to $\phi\,\circ$ , \emph{composition with} $\phi$
(\cite{Sanhan}). The two
tensor-extended compositions are consistent.

\subsection{Quantum stochastics
(\cite{Partha}{\rm ,} \cite{Meyer}{\rm ; we follow}
\cite{Lgreifswald}{\rm ,} \cite{LSqsde}{\rm , modified for abstract
spaces})}
\emph{Fix now, and for the rest of the paper}, a complex Hilbert space
$\kil$ which we refer to as the \emph{noise dimension space},
and let $\kilhat$ denote the orthogonal sum $\Comp\oplus\kil$. Whenever $c\in \kil$,
$\chat:=\binom{1}{c}\in \kilhat$; for $E\subset\kil$,
$\widehat{E}:= \Lin\{\chat:c \in E\}$ and when $g$ is a function with
values in $\kil$, $\ghat$ denotes the corresponding function with values
in $\kilhat$, defined by $\ghat(s):= \widehat{g(s)}$.
Let $\FFock_I$ denote the symmetric Fock space over $L^2(I;\kil)$,
dropping the subscript when the interval $I$ is all of $\Rplus$.
For any dense subspace $D$ of $\kil$ let
$\Step_D$ denote the linear span of $\{d_{[0,t[}: d\in D, t\in\Rplus\}$
in $L^2(\Rplus;\kil)$ (we always take these right-continuous versions) and
let $\ExpsD$ denote the linear span of $\{\ve(g): g\in\Step_D\}$ in
$\FFock$, where $\ve(g)$ denotes the exponential vector
$\big((n!)^{-\frac{1}{2}}g^{\ot n}\big)_{n\geq 0}$.
The subscript $D$ is dropped when $D=\kil$.
We usually drop the tensor symbol and denote
simple tensors such as $v\ot\ve(f)$ by $v\ve(f)$.
Also define
\begin{equation} \label{Delta QS}
e_0 := \binom{1}{0}\in\kilhat \text{ and }
\dQS := P_{\{0\}\oplus\kil}
=
\begin{bmatrix}
0 &  \\  & I_{\kil}
\end{bmatrix}
\in B(\kilhat).
\end{equation}

The basic objects we consider in this paper are completely bounded quantum
stochastic mapping processes on operator spaces. These are time-indexed
families of completely bounded maps $\{k_t: t \geq 0\}$ from an operator
space to the algebra of bounded operators on $\init \ot \FFock$, for a
Hilbert space $\init$, satisfying standard
adaptedness and measurability conditions. For technical reasons we also
need to consider mapping processes whose values are (at least, a priori)
unbounded operators.
The crucial point here is
that the naturally arising operators have `bounded slices':
for any vectors $\ve, \ve' \in \Exps$ the maps
\[
v \to
\big(
I_{\init} \ot \la \ve' | \big)
\, k_t(v) \,
\big( I_{\init} \ot | \ve \ra \big)
\]
($t\in\Rplus$) have values in $B(\init)$, and are (completely) bounded,
even though the global maps $k_t$ may not be -- more precisely they have
(completely) bounded columns (see Property 2, following
Theorem~\ref{thm: existence}).
This point of view, where each $k_t$ is taken to be a family of maps
indexed by pairs of exponential vectors, allows the replacement of
$B(\init)$ by an abstract operator space and, once the somewhat technical
definitions below are accepted, leads to a development of the theory
which is straight-forward and effective with more transparent proofs.
This said, to follow the arguments it is safe to keep in mind sesquilinear
maps induced by mapping processes in the familiar sense.

Let $\Vil$ and $\Wil$ be operator spaces. In this paper we denote by
$\Proc (\Vil \Pto \Wil)$ the collection of families $k= (k_t)_{t\geq 0}$
of maps in
\[
L\big(\Vil ; L\big(\Exps ; CB ( \la \FFock | ; \Wil ) \big)\big)
\subset SL \big( \Exps , \Exps ; L(\Vil ; \Wil ) \big)
\]
satisfying the following \emph{measurability} and \emph{adaptedness}
conditions
\begin{align*}
& s \mapsto k^{\ve',\ve}_s \text{ is pointwise weakly measurable, and } \\
& k^{\ve',\ve}_t = \la \ve'_2, \ve_2 \ra \, k^{\ve'_1, \ve_1}_t,
\end{align*}
for $\ve = \ve(f) , \ve' = \ve (f') \in \Exps$ and $t \in \Rplus $, where
$\ve_1 = \ve(f_{[0,t[})$
and $\ve_2 = \ve (f_{[t,\infty[})$ with $\ve'_1$ and $\ve'_2$ defined in
the same way for $f'$.
When $\Wil = \Comp$ (as is the case for quantum stochastic convolution
cocycles) we write $\cProc(\Vil)$ instead of $\Proc (\Vil \Pto \Comp)$.
Then $k_t \in L \big( \Vil ; \Op (\Exps) \big)$ for each $t\geq 0$
and, in terms of the exponential property of Fock space:
$\FFock = \FFock_{[0,t]}\ot\FFock_{[t,\infty[}$, adaptedness reads
\[
k_t (x) \ve (f) = u_t \ot \ve (f_{[t,\infty[}) \text{ where }
u_t = k_t (x) \ve (f_{[0,t[}) \in \FFock_{[0,t]}.
\]
Here the following Banach space identifications are used:
\[
CB\left(\la\FFock |;\Comp\right) = B\left(\la\FFock |;\Comp\right)
=\FFock.
\]
When $k_t$ is viewed as a map in
$L\big(\Vil,\Exps; CB(\la\FFock |;\Wil )\big)$
we use the notation $k_{t,|\ve\ra}(x)$. Note that if
$k\in\Proc(\Vil\Pto\Yil)$ for a concrete operator space $\Yil$ then,
invoking the complete isometry~\eqref{orawan},
$k_{t,|\ve\ra}\in L\big(\Vil;\Yil\otM |\FFock\ra\big)$.

Processes $k$ and $j$ are identified if, for all
$\ve', \ve \in \Exps, x \in \Vil$ and $\varphi \in \Wil^*$,
the scalar-valued functions
$t \mapsto \varphi \circ k_t^{\ve',\ve} (x)$ and $t \mapsto \varphi \circ j_t^{\ve', \ve} (x)$
agree almost everywhere. We also denote by $\Procdagger (\Vil \Pto \Wil)$
the subspace of processes $k$ for which
\[
\text{ each map }
\la \ve' | \mapsto k^{\ve, \ve'}_t (x)
\text{ is completely bounded } \la \Exps| \to \Wil
\]
($\ve \in \Exps, x \in \Vil ,t \in \Rplus$).
Then, for $k \in \Procdagger (\Vil \Pto \Wil)$,
\[
(k^\dagger )^{\ve',\ve}_t := (k^{\ve, \ve'}_t )^\dagger
\]
defines a process
$k^\dagger \in \Procdagger (\Vil^\dagger \Pto \Wil^\dagger )$, where
$\dagger$ denotes
conjugate operator space. When the operator space $\Wil$ is concrete this
amounts to the usual notion of adjoint(able) process.
\emph{Complete boundedness} for a process $k \in \Proc (\Vil \Pto \Wil)$
means
\[
k_t \in CB \big( \la \FFock | , | \FFock \ra; CB(\Vil ; \Wil ) \big)
\subset L \big( \Vil , \Exps ; CB( \la \FFock | ; \Wil ) \big)
\]
for each $t \in \Rplus$. Thus $\Proccb (\Vil \Pto \Wil)$, the class of
such processes, is a subspace of
$\Procdagger (\Vil \Pto \Wil)$. The natural inclusion
\begin{equation} \label{worth noting 1}
CB \big( \Vil ; \Wil \ot B(\FFock)\big) \subset
CB \big( \la \FFock | , | \FFock \ra ; CB(\Vil ; \Wil ) \big)
\end{equation}
and, for a concrete operator space $\Yil$, the natural identification
\begin{equation} \label{worth noting 2}
CB \big( \la \FFock | , | \FFock \ra ; CB(\Vil ; \Yil ) \big)
= CB\big( \Vil ; \Yil \otM B(\FFock) \big)
\end{equation}
(\cite{Sanhan}) are both worth noting here (cf.\,~\eqref{sp incl cb});
they explain the terminology.

We need two further properties for processes: $k\in\Proc (\Vil\Pto\Wil)$
is \emph{weakly initial space bounded} if
\[
k^{\ve',\ve}_t : \Vil \to \Wil \text{ is bounded }
\]
$(\ve , \ve' \in \Exps, t \in \Rplus)$ and is \emph{weakly regular} if
further
\[
\sup \big\{ \| k_s^{\ve', \ve} \| : 0 \leq s \leq t \big\} < \infty ,
\]
for all $t\geq 0$.
We shall be dealing with quantum stochastic differential equations of the form
\begin{equation} \label{QSDE}
dk_t = k_t \fullcomp d\Lambda_\phi (t), \quad
k_0 = \iota_{\FFock} \circ \kappa ,
\end{equation}
where
$\phi \in CB\big( \Vil ; \Vil \ot B(\kilhat) \big)$
and
$\kappa\in CB(\Vil ;\Wil )$.
Here the natural inclusion
\begin{equation} \label{both relevant 1}
CB \big( \Vil ; \Vil \ot B(\kilhat ) \big) \subset
CB \big( \la \kilhat | , | \kilhat \ra ; CB (\Vil ) \big)
\end{equation}
and, for a concrete operator space $\Yil$ in $B(\init)$, the natural
complete isometries
\begin{equation} \label{both relevant 2}
CB \big( \la \kilhat | , | \kilhat \ra ; CB(\Yil) \big)
= CB\big( \Yil; \Yil \otM B(\kilhat) \big)
\subset
CB\big( \Vil ; B(\init \ot \kilhat) \big)
\end{equation}
are relevant
(cf.\,~\eqref{worth noting 1} and~\eqref{worth noting 2}).
A process $k \in \Proc (\Vil \Pto \Wil )$ is a
\emph{weak}
solution of the QS differential equation~\eqref{QSDE} if
\begin{align*}
& s \mapsto k^{\ve', \ve}_s \circ \phi^{\zeta',\zeta} (x)
   \text{ is weakly continuous, and } \\
& k^{\ve', \ve}_t (x) = \la \ve', \ve \ra \kappa (x)
  + \text{w-}\!\!\int^t_0 k^{\ve', \ve}_s \circ \phi^{\wh{f}' (s) , \fhat (s)} (x) \, ds
\end{align*}
$(\zeta,\zeta'\in\kilhat,
\ve = \ve (f), \ve' = \ve (f') \in \Exps , x \in \Vil ,t \in \Rplus )$; it is so-called
due to the First Fundamental Formula of quantum stochastic calculus.
The theorem we need
is the following special case of Proposition 3.5 in~\cite{LSqsde} which
generalises~\cite{LWfeller}
to allow nontrivial initial conditions and abstract spaces.

\begin{thm}\label{thm: existence}
Let $\phi \in CB\big( \Vil ; \Vil \ot B(\kilhat) \big)$
and $\kappa \in CB(\Vil ;\Wil)$. Then there is a unique
weakly regular weak solution of the quantum stochastic
differential equation~\eqref{QSDE}.
\end{thm}

\begin{notn}
$k^{\kappa, \phi}$, simplifying to $k^\phi$ for the case where
$\Vil =\Wil$ and $\kappa = \id_{\Wil}$.
\end{notn}

We list key properties of the solution processes needed in this paper
next (see~\cite{LSqsde}).
Let $k = k^{\kappa , \phi}$ for $\kappa$ and $\phi$ as above.
\begin{enumerate}[1.]
\item
$k \in \Procdagger (\Vil \Pto \Wil)$ and $k^\dagger = k^{\kappa^\dagger,
\phi^\dagger}$.
\item \label{stHol1}
For all $\ve \in \Exps$ and $t \in \Rplus$,
$k_{t, | \ve \ra} \in CB \big( \Vil ; CB ( \la \FFock | ; \Wil ) \big)
=
CB \big(\la \FFock |; CB(\Vil;\Wil)\big)$
(the process has completely bounded columns)
and the map $s \mapsto k_{s, | \ve \ra}$ is locally H\"{o}lder-continuous with
exponent $1/2$. Moreover if
$\phi (\Vil)\subset\Vil\ot B(\kilhat )$,  then $k$ satisfies
\[
k_{t,|\ve\ra}(\Vil )\subset\Wil\ot |\FFock\ra
\]
(in terms of the inclusion~\eqref{abstract model}).
\item
If $\kappa = \kappa_2 \circ \kappa_1$ where
$\kappa_1 \in CB(\Vil ; \Uil)$ and
$\kappa_2 \in CB(\Uil ; \Wil)$ then
\[
k_{t,| \ve \ra} =
\kappa_2 \fullcomp \wt{k}_{t,| \ve \ra}
:x\mapsto \kappa_2 \circ \wt{k}_{t,| \ve \ra}(x)
\text{ where }
\wt{k} = k^{\kappa_1 , \phi}
\]
$(t \in \Rplus , \ve \in \Exps )$. If the process $\wt{k}$ is completely bounded then so is $k$
and we have the identity
\[
k_t = \kappa_2 \fullcomp k^{\kappa_1 ,\phi}_t
\]
($t \in \Rplus$). In particular,
\[
k_{t,| \ve \ra} = \kappa \fullcomp k^\phi_{t,| \ve \ra}
\quad
(\text{resp. } k_t = \kappa \fullcomp k^\phi_t
\text{ when } k^\phi \in \Proccb (\Vil \Pto \Wil ) \big) .
\]
\item
The following useful `form representation' holds
\begin{equation} \label{form repn}
k^{\ve', \ve}_t(x) =
\la \ve', \ve \ra \ \text{w-}\!\!\int_{\Gamma_{[0,t]}} \,
d \sigma \ \Omega_\sigma \circ \tau_{\# \sigma} \circ \kappa \fullcomp \phi^{\fullcomp \# \sigma}(x)
\end{equation}
(weak integral) where $\Omega_\sigma := \omega_{\xi' , \xi}$ for $\xi = \pi_{\fhat} (\sigma)$
and $\xi' = \pi_{\wh{f'}} (\sigma)$, when
$\ve = \ve (f)$ and $\ve' = \ve (f')$ and, for $n\in\Nat$, $\tau_n$
is the permutation of tensor components which reverses the order of the $\kilhat$'s.

\item
Let $\ve , \ve' \in \Exps$ and $t \in \Rplus$. If each
$k_{t, | \ve \ra}$ is $\Wil \ot | \FFock \ra$-valued
then
\[
k^{\ve', \ve}_t =
(\id_{\Wil} \ot \lambda_{\ve'} ) \circ k_{t, | \ve \ra },
\]
and if further $k \in \Proccb (\Vil \Pto \Wil)$ and each $k_t$ maps
$\Vil$ into the spatial tensor product $\Wil \ot B(\FFock)$, then
\[
k_{t,| \ve \ra} = ( \id_{\Wil} \ot \rho_{\ve} ) \circ k_t.
\]
Here $\lambda_{\ve'} : | \FFock \ra \to \Comp$ denotes left
multiplication by $\la \ve' |$, and
$\rho_{\ve} : B(\FFock) \to | \FFock \ra$ right multiplication by $| \ve \ra$.

\item
In fact $k$ is a \emph{strong solution} of the QS differential equation,
meaning that the integral equation
\[
k_t = \iota_{\FFock} \circ \kappa + \int^t_0 k_s \fullcomp d \Lambda_\phi
(s)
\]
is valid in a `strong sense'.

\item
The process $k$ is expressible in
terms of the multiple QS integral operation:
\[
k_t = \Lambda_t \circ \upsilon \text{ where }
\upsilon = \big( \upsilon_n \big)_{n\geq 0} \text{ and }
\upsilon_n = \tau_n \circ \kappa \fullcomp \phi^{\fullcomp n}
\]
(cf.\ Property 4).
\item
When $\Vil=\Wil$ and $\kappa = \id_{\Wil}$, $k_0 = \iota_{\FFock}$ and $k$
enjoys the following weak cocycle property: for $s,t\in\Rplus$,
$\ve = \ve (f)$ and  $\ve' = \ve(f')$ in $\Exps$,
\[
k_{s+t}^{\ve',\ve}
=
\la \ve'_3 , \ve_3 \ra \,
k_s^{\ve'_1, \ve_1}  \circ
k_t^{\ve'_2,\ve_2},
\]
where
\begin{equation} \label{epsilons}
\ve_1 = \ve (f_{[0,s[}) , \
\ve_2 = \ve (S_s^* f_{[s,s+t[}) \text{ and }
\ve_3 = \ve (f_{[s+t,\infty[}) ,
\end{equation}
$\ve'_1$, $\ve'_2$ and $\ve'_3$ being
defined similarly with $f'$ in place of $f$,
$(S_t)_{t\geq0}$ is the one-parameter semigroup of right shifts on
$L^2(\Rplus;\kil)$ and $(\sigma_t)_{t\geq 0}$ is the induced endomorphism
semigroup on $B(\FFock)$, ampliated to $SL\left(\Exps,\Exps;\Wil\right)$.
The cocycle property simpifies to
\[
k_{s+t} = k_s \fullcomp \sigma_s \circ k_t
\]
when $k$ is completely bounded.
\end{enumerate}
Property 8, namely the fact that processes of the form $k^\phi$ are weak QS
cocycles (also called Markovian cocycles), has a converse --- subject to
certain constraints, weak QS cocycles are necessarily of this form. The main
results in this direction concern
completely positive, contractive QS cocycles
on a unital $C^*$-algebra and are collected next ---
they originate in \cite{grandmother};
a direct proof is given in~\cite{LWdirect}.
\begin{thm}
[{\rm [}$\text{LW}_{\!1-3}${\rm ]}]
\label{CPstandard}
Let $\alg$ be a unital $C^*$-algebra
and let $k\in \Proc(\alg \Pto \alg)$.
Then the following are equivalent\tu{:}
\begin{rlist}
\item
$k$ is a Markov-regular, completely positive and contractive QS
cocycle on  $\alg$\tu{;}
\item
$k=k^{\phi}$ where
$\phi \in CB \big( \la\kilhat |,|\kilhat\ra;CB(\alg)\big) = CB \big(\alg; CB(\la\kilhat |,|\kilhat\ra; \alg) \big)$ satisfies
$\phi (1) \leq 0$
and, in any faithful, nondegenerate representation,
$\phi$ may be decomposed as follows\tu{:}
\begin{equation}   \label{CPC0}
 \phi(x) =
\Psi (x) -  x \ot \Projk - \big(x \ot | e_0 \ra\big)J -
J^* \big( x \ot \la e_0 | \big)
\end{equation}
\tu{(}$x\in\alg$\tu{)}, for some map
$\Psi\in CP\big(\alg ;\alg''\olot B(\kilhat )\big)$ and operator
$J \in \alg'' \olot \la \kilhat |$.
\end{rlist}
\end{thm}
\noindent
Here, with respect to the
representation, $\alg''$ denotes the double commutant of $\alg$, $\olot$ denotes the
ultraweak tensor product and $\phi(x) \in \alg \otM B(\kilhat)$.
The form~\eqref{CPC0} taken by the stochastic generator (see also
\cite{Slava}) generalises the
Christensen-Evans Theorem on the generators of norm continuous completely
positive contractive semigroups (\cite{chrisevans}).

\section{Operator space coalgebras} \label{sec: B}

In this section we adapt the basic notions of coalgebra to the category of operator spaces,
and consider convolution semigroups of functionals in this context.
Three structures are considered: operator space coalgebras, operator
system coalgebras and $C^*$-bialgebras, corresponding to the three levels
of question addressed in this paper, namely linear, positivity-preserving
and algebraic. In fact a hybrid structure, called $C^*$-hyperbialgebra,
plays a more prominent role than operator system coalgebras do.

\begin{defn}
An \emph{operator space coalgebra} is an operator space $\Cil$
equipped with complete contractions
$\counit : \Cil\to\Comp$ and $\Com : \Cil \to\Cil \ot \Cil$,
called the \emph{counit} and \emph{coproduct} respectively, satisfying
\begin{enumerate}[\,\,\,\,(OSC1)]
\item  \label{coassoc}
$(\Com \ot \id_{\Cil} ) \circ \Delta =
(\id_{\Cil} \ot \Com) \circ \Com  $  \quad
(coassociativity),
\item   \label{counit prop}
$(\counit \ot \id_{\Cil} ) \circ \Delta = \id_{\Cil} = (\id_{\Cil} \ot \counit ) \circ \Com$
\quad (counit property);
\end{enumerate}
it is an \emph{operator system coalgebra} if $\Cil$ is an operator system
and
\begin{enumerate}[\,\,\,\,(OSyC)]
\item
$\counit$ and $\Com$ are both unital and completely positive;
\end{enumerate}
a \emph{$C^*$-hyperbialgebra} if $\Cil$ is a unital $C^*$-algebra and
\begin{enumerate}[\,\,\,\,($C^*$-Hy)]
\item
$\counit$ is a character (i.e. it is nonzero and multiplicative) and
$\Com$ is unital and completely positive;
\end{enumerate}
and finally it is a \emph{$C^*$-bialgebra} if $\Cil$ is a (unital)
$C^*$-algebra and
\begin{enumerate}
[($C^*$-Bi)]
\item
$\counit$ and $\Com$ are both unital and *-homomorphic.
\end{enumerate}
\end{defn}

By (OSC\ref{coassoc}), $\Com^{\fullcomp 2}$ is defined unambiguously,
as is $\Com^{\fullcomp n}$ for all $n \in \Nat$, and we define
$\Com^{\fullcomp 0} := \id_{\Cil}$.
Similarly (OSC\ref{counit prop}) gives $\counit \fullcomp \Com =
\id_{\Cil}$, with
unambiguous meaning. An operator space coalgebra is \emph{cocommutative} if
$\Com = \Com^{\opposite}$ where
\[
\Com^{\opposite} := \tau \circ \Com,
\]
$\tau$ being the tensor flip on $\Cil \ot \Cil$. The \emph{opposite}
operator space coalgebra results from replacing $\Com$ by
$\Com^{\opposite}$.

An operator space coalgebra is thus typically \emph{not} a coalgebra in
the algebraic sense (\cite{Sweedler}) since the coproduct is not required
to map $\Cil$ into $\Cil \odot\Cil$.
A (unital) $C^*$-bialgebra is a $C^*$-hyperbialgebra $\Bil$
whose coproduct is also multiplicative (and thus a unital *-homomorphism).
Some authors (for example \cite{vaes}) drop the unital condition on
$C^*$-bialgebras,
and require instead the counit to be a nondegenerate *-homomorphism into
the multiplier algebra $M(\Cil \ot \Cil)$.
The asymmetry in the definition of $C^*$-hyperbialgebra --- whereby
$\counit$ is required to be multiplicative but $\Com$ only
to be completely positive --- is motivated by the example of compact
quantum hypergroups (\cite{cqh}). Multiplicativity of the counit is used
extensively in characterising generators of completely positive convolution
cocycles (Section~\ref{sec: F} and \cite{Sthesis}). Finally note that the
conjugate operator space of an operator space coalgebra has natural
operator space coalgebra structure.

\subsection{Convolution}
For an operator space coalgebra $\Cil$ and operators spaces $\Vil_1$ and $\Vil_2$,
the \emph{convolution} of $\varphi_1 \in CB(\Cil ; \Vil_1)$ and $\varphi_2 \in CB(\Cil ;\Vil_2)$
is defined by
\[
\varphi_1 \fullstar \varphi_2 := (\varphi_1 \ot \varphi_2) \circ \Com \in
CB(\Cil ; \Vil_1 \ot \Vil_2 ).
\]
It is easily seen that convolution is associative
(in the same sense as the spatial tensor product is)
and enjoys submultiplicativity and unital properties:
\begin{align}
&(\varphi_1 \fullstar \varphi_2) \fullstar \varphi_3 =
\varphi_1 \fullstar (\varphi_2 \fullstar \varphi_3)
\label{star associativity} \\
& \| \varphi_1 \fullstar \varphi_2 \|_{\cb} \leq
\| \varphi_1 \|_{\cb} \| \varphi_2 \|_{\cb}, \text{ and } \notag \\
& \counit \fullstar \varphi = \varphi = \varphi \fullstar \counit . \label{conv unit}
\end{align}
In particular, $(\Cil^* , \fullstar)$ is a unital Banach algebra.
For $n \in \Nat$ and $\varphi_1, \ldots ,\varphi_n \in CB(\Cil ; \Vil)$,
$n$-fold convolution is
defined via $n$-fold tensor products:
\begin{equation} \label{phi star n}
\varphi_1\fullstar \cdots \fullstar \varphi_n
= (\varphi_1 \ot \cdots \ot \varphi_n) \circ \Com^{\fullcomp (n-1)} .
\end{equation}
We also define
$\varphi^{\fullstar 0} := \counit$, which is consistent with~\eqref{conv unit}.

Given an operator space coalgebra $\Cil$, each operator space $\Vil$
determines maps
\begin{align*}
&R_{\Vil} : CB(\Cil;\Vil) \to CB(\Cil ;\Cil \ot \Vil),
     \quad \varphi \mapsto (\id_{\Cil} \ot \varphi) \circ \Com; \\
&E_{\Vil} :CB(\Cil;\Cil \ot \Vil) \to CB(\Cil ;\Vil),
     \quad \phi \mapsto (\counit \ot \id_{\Vil} )\circ \phi .
\end{align*}
Thus the action of $R_{\Vil}$ is convolve with the identity map on $\Cil$, putting the
argument on the \emph{right}, and that of $E_{\Vil}$ is compose in the
tensor-extended sense with the counit:
\[
R_{\Vil} \varphi = \id_{\Cil} \fullstar \varphi, \text{ and } E_{\Vil} \phi = \counit \fullcomp \phi .
\]
In the noncocommutative case we are therefore making a choice here.
We abbreviate $R_{\Comp}$ to $\cR$.

The basic properties of these maps are collected below. They are all
easily proved
from the definitions, noting that under the completely isometric
identification $\Mat_n (\Cil \ot \Vil)= \Cil \ot \Mat_n (\Vil)$,
\[
(R_{\Vil})^{(n)} = R_{\Mat_n (\Vil)} .
\]

\begin{propn} \label{thm: R and E}
Let $\Cil$ be an operator space coalgebra, and let $\Vil_1$, $\Vil_2$ and
$\Vil$ be operator spaces.
\begin{alist}
\item
$R_{\Vil}$ and $E_{\Vil}$ are complete isometries satisfying
\[
E_{\Vil} \circ R_{\Vil} = \id_{CB(\Cil ;\Vil)} .
\]
\item
If $\varphi_1 \in CB(\Cil ;\Vil_1)$ and $\varphi_2 \in CB(\Cil ; \Vil_2)$
then
\[
R_{\Vil_1 \ot \Vil_2} (\varphi_1 \fullstar \varphi_2)
= R_{\Vil_1} \varphi_1 \fullcomp R_{\Vil_2} \varphi_2 .
\]
\item
If $\varphi\in CB(\Cil;\Vil)$ then
\[
R_{\Vil^\dagger} (\varphi^\dagger) = (R_{\Vil} \varphi)^\dagger .
\]
\end{alist}
\end{propn}

\begin{rem}
Noting that $\counit =E_{\Comp}(\id_{\Cil})$ and $\Com = R_{\Cil}
(\id_{\Cil})$, it is clear
that operator space coalgebras could be axiomatised in terms of $R$- and
$E$-maps in lieu of $\Com $ and $\counit$.
\end{rem}

Write $\CBcoproduct (\Cil ;\Cil \ot \Vil)$ for $\Ran R_{\Vil}$.

\begin{cor} \label{C}
For each operator space $\Vil$, $R_{\Vil}$ determines a complete isometry
of operator spaces
\[
CB(\Cil ; \Vil) \cong \CBcoproduct (\Cil ; \Cil \ot \Vil ),
\]
by corestriction. In case $\Vil = \Comp$ this gives
an isometric isomorphism of unital Banach algebras
\[
(\Cil^* , \fullstar) \cong \big( \CBcoproduct (\Cil ), \circ \big) .
\]
\end{cor}

A further noteworthy consequence is the following identity.

\begin{cor} \label{C+}
In $\CBcoproduct (\Cil ; \Cil \ot \Mat_n )$,
\[
\| \phi\|_{\cb} = \| \phi^{(n)} \| .
\]
\end{cor}

\begin{proof}
Let $\phi \in \CBcoproduct (\Cil ; \Cil \ot \Mat_n )$, say $\phi = R_{\Mat_n} \varphi$.
Then $\varphi\in CB(\Cil; \Mat_n)$ so
\[
\| \phi \|_{\cb} = \| \varphi \|_{\cb} = \| \varphi^{(n)} \| = \| \counit \fullcomp \phi^{(n)} \|
\leq \| \phi^{(n)} \|.
\]
The result follows.
\end{proof}
In particular, in $\CBcoproduct (\Cil)$ the completely bounded norm coincides with the bounded
operator norm. As a result $\CBcoproduct (\Cil)$ is a closed subspace of $B(\Cil)$.
The next proposition collects the structure-preserving properties of
$R_{\Vil}$ and its inverse, under a number of pertinent assumptions on
$\Cil$ and $\Vil$.

\begin{propn} \label{Requiv}
Let $\Cil$ be an operator space coalgebra and $\Vil$ an operator space, let
$\varphi \in CB(\Cil ;\Vil)$ and
$\phi = R_{\Vil}\varphi \in \CBcoproduct (\Cil ; \Cil \ot \Vil
)$.
\begin{alist}
\item
The map $\phi$ is completely contractive if and only if $\varphi$ is.
\item
If $\Cil$ is an operator system coalgebra and $\Vil$ is an operator system
then $\phi$ is real \tu{(}respectively, completely positive, or unital\tu{)} if and
only if $\varphi$ is.
\item
If $\Cil$ is a $C^*$-bialgebra and $\Vil$ is a $C^*$-algebra then $\phi$
is multiplicative if and only if $\varphi$ is.
\end{alist}
\end{propn}

A \emph{convolution semigroup of functionals} on an operator space coalgebra
$\Cil$ is a one-parameter family $\lambda =(\lambda_t)_{t \geq 0}$ in
$\Cil^*$ satisfying
\[
\lambda_0 = \counit \text{ and } \lambda_{s+t} = \lambda_s \fullstar \lambda_t .
\]
In other words a convolution semigroup of functionals on $\Cil$ is a
one-parameter semigroup in the unital
algebra $\big(\Cil^*,\fullstar\big)$.

\begin{propn}
Let $\Cil$ be an operator space coalgebra. The map
$\lambda \mapsto P := (\cR \lambda_t)_{t \geq 0}$
is a bijection from the set of convolution semigroups of functionals
on $\Cil$ to the set of one-parameter semigroups in
$\CBcoproduct (\Cil)$. Moreover, the conditions in (a) below are equivalent, and so are the conditions in (b)\tu{:}
\begin{alist}
\item   \label{part a}
\begin{alist}
\item
\label{lim}
$\lambda_t  \to \counit$ pointwise as $t \to 0$\tu{;}
\item
\label{c0 semigroup}
$P$ is a $C_0$-semigroup on $\Cil$.
\end{alist}
\item
\label{part b}
\begin{rlist}
\item
$\lambda$ is norm continuous in $t$\tu{;}
\item
$P$ is  norm continuous in $t$\tu{;}
\item
$P$ is $\cb$-norm continuous in $t$;
\item $P$ has a completely bounded generator.
\end{rlist}
\end{alist}
\end{propn}

\begin{proof}
The first part follows from Corollary~\ref{C}.

Since $\counit \circ P_t = \lambda_t$,
(\ref{c0 semigroup}) implies (\ref{lim}). Suppose
therefore that (\ref{lim}) holds. Then, for any $\varphi \in \Cil^*$,
\[
\varphi \circ P_t = \lambda_t \circ (\varphi \ot \id_{\Cil}) \circ \Com
\text{ and } \counit \circ (\varphi \ot \id_{\Cil} ) \circ \Com = \varphi ,
\]
so $P_t x \to x$ weakly as $t \searrow 0$, for all $x \in \Cil$. But this implies that $P$ is
strongly continuous (\cite{Davies}, Proposition 1.23) and thus a
$C_0$-semigroup, so (\ref{c0 semigroup}) holds.

By Corollary~\ref{C+}
\[
\| P_t -\id_{\Cil} \|_{\cb} = \| \lambda_t - \counit \| = \| P_t -
\id_{\Cil} \| ,
\]
and so (\ref{part b}) follows.
\end{proof}
Thus each norm-continuous convolution semigroup of functionals $\lambda$
on $\Cil$ has a \emph{generator}:
\[
\gamma := \lim_{t \searrow 0} t^{-1} (\lambda_t - \counit )
\]
from which the convolution semigroup of functionals may be recovered
\[
\lambda_t =
\exp_{\fullstar} t \gamma :=
\sum_{n \geq 0} \frac{t^n}{n!} \gamma^{\fullstar n} .
\]
The corresponding one-parameter semigroup on $\Cil$ has completely bounded
generator:
\[
\cR\lambda_t = e^{t\tau}, \text{ where }
\tau = \cR\gamma \in \CBcoproduct (\Cil) .
\]

\section{Operator space coalgebraic QS differential equations} \label{sec: D}

In this section we consider operator space coalgebraic quantum stochastic differential
equations with completely bounded coefficients, and relate their solutions to those of
standard QS differential equations by means of $R$-maps. In particular we show that the complete boundedness
property is preserved when moving between these two kinds of solutions.
\emph{For this section
$\Cil$ is a fixed operator space coalgebra.}

Let $\varphi \in CB \big( \Cil ; B(\kilhat)\big)$. A weakly initial space
bounded process $k \in \cProc (\Cil)$ is a \emph{weak solution}
of the operator
space coalgebraic
quantum stochastic differential equation
\begin{equation} \label{OScoalgQSDE}
dk_t = k_t \fullstar d \Lambda_\varphi (t) , \quad
k_0 = \iota_{\FFock} \circ \counit ,
\end{equation}
if
\begin{equation} \label{w2coalg}
s \mapsto\big(k^{\ve',\ve}_s \fullstar \varphi^{\zeta',\zeta}\big)(x)
\text{ is continuous, and }
\end{equation}
\begin{equation} \label{w3coalg}
k^{\ve',\ve}_t (x) = \la \ve', \ve \ra \kappa (x) +
  \int^t_0 (k^{\ve',\ve}_s \fullstar \varphi^{\wh{f'} (s), \fhat (s)} ) (x) \, dx
\end{equation}
for $\zeta,\zeta'\in\kilhat, \ve=\ve(f),\ve'=\ve(f')\in\Exps, x\in\Cil$
and $t \in \Rplus$.

\begin{rem}
By the Banach-Steinhaus Theorem
\[
\sup \Big\{ \big| (\omega_{\ve (f'), \ve(f)} \circ k_s) \fullstar
(\omega_{\wh{f'}(s),\wh{f}(s)} \circ \varphi ) (x) \big| \,
\Big| \, x \in \Cil, \| x\| \leq 1, s \in [0,t] \Big\} < \infty
\]
for each $f,f'\in\Step$ and $t\in\Rplus$.
It follows therefore that weak solutions of the operator space
coalgebra QS differential equation in the above sense are automatically
weakly regular.
\end{rem}

\begin{thm} \label{thm: coalg existence}
Let $\varphi \in CB\big(\Cil ;B(\kilhat) \big)$. Then the operator space
coalgebraic quantum stochastic differential equation~\eqref{OScoalgQSDE}
has a unique weak solution.
\end{thm}

\begin{proof}
Let $k \in \cProc ( \Cil )$ be weakly regular. Then
\[
k^{\ve',\ve}_t \fullstar \varphi^{\zeta',\zeta} =
k^{\ve',\ve}_t \circ \phi^{\zeta',\zeta}
\]
$(\ve , \ve' \in \Exps , \zeta , \zeta' \in \kilhat , t \in \Rplus )$
where $\phi \in R_{B(\kilhat)} \varphi$.
It follows that $k$  weakly satisfies the operator space coalgebraic
QS differential equation~\eqref{OScoalgQSDE}
if and only if $k$ weakly satisfies the operator space QS differential
equation
$dk_t =k_t \fullcomp d\Lambda_\phi (t) , \,
k_0 = \iota_{\FFock} \circ \counit$. Since
$\phi \in CB\big(\Cil ;\Cil \ot B(\kilhat) \big)$
and $\counit \in \Cil^* = CB(\Cil ; \Comp)$ the result therefore follows
from Theorem~\ref{thm: existence},
and the automatic weak regularity of weak solutions of~\eqref{OScoalgQSDE}.
\end{proof}

\begin{notn}
We denote the unique weak solution
of~\eqref{OScoalgQSDE}, for completely bounded $\varphi$, by $l^\varphi$.
From the above proof we see that $l^\varphi = k^{\counit , \phi}$         where
\[
\phi = R_{B (\kilhat)} \varphi \in CB \big( \Cil ; \Cil \ot B(\kilhat ) \big) .
\]
\end{notn}

Note that Proposition~\ref{thm: R and E} implies that
\[
\counit \fullcomp \phi^{\fullcomp n} =
\counit \fullcomp R_{B(\kilhat)^{\ot n}} \varphi^{\fullstar n} =
\varphi^{\fullstar n}, \quad n \geq 0 .
\]

The properties of solutions of operator space QS differential equations
listed in Section~\ref{sec: A} entail the following for $l = l^\varphi$
where $\varphi \in CB\big(\Cil ; B(\kilhat) \big)$:

\begin{enumerate}[1$'$.]
\item
$l \in \cProcdagger ( \Cil )$ and $l^\dagger = l^\psi$ where
$\psi = \varphi^\dagger \in CB \big( \Cil^\dagger ; B(\kilhat ) \big)$.
\item \label{stHol2}
$l_{t ,| \ve \ra} \in CB \big( \Cil ;| \FFock \ra \big)$ and the map $s \mapsto l_{s, | \ve \ra}$
is locally H\"{o}lder continuous with exponent $\frac12$, moreover
$k^\phi_{t, | \ve \ra}(\Cil)\subset\Cil\ot|\FFock\ra$ for all $ \ve \in \Exps$ and $t \in \Rplus$.
\item
Since $l = k^{\counit , \phi}$, where $\phi = R_{B(\kilhat)} \varphi$,
\begin{equation} \label{l phi epsilon}
l_{t, | \ve \ra} = \counit \fullcomp k^\phi_{t, | \ve \ra}
\end{equation}
$( \ve \in \Exps , t \in \Rplus )$; also
if $k^\phi$ is completely bounded then $l$ is too and
\[
l_t = \counit \fullcomp k_t^\phi, \quad t \in \Rplus.
\]
\item
In the notation of Property 4,
\begin{equation} \label{again from repn}
l^{\ve',\ve}_t = \la \ve',\ve \ra \int_{\Gamma_{[0,t]}}
d\sigma \ \Omega_\sigma \circ \tau_{\# \sigma} \circ \varphi^{\fullstar \# \sigma}
\end{equation}
for $\ve = \ve(f) , \ve' \in \ve (f') \in \Exps$ and $t \in \Rplus$.
\item
In the notation of Property 5,
\[
l^{\ve',\ve}_t = \lambda_{\ve'} \circ l_{t, | \ve \ra}
\]
and, if $l$ is completely bounded so that $l_t \in CB\big(\Cil ; B(\FFock ) \big)$, then
\[
l_{t, | \ve \ra} = \rho_\ve \circ l_t
\]
$( \ve , \ve' \in \Exps , t \in \Rplus )$.
\item
$l$ is a strong solution of the operator space coalgebraic QS
differential equation:
\[
l_t =
 \iota_{\FFock} \circ \counit +
\int^t_0 l_s \fullstar \, d \Lambda_\varphi(s)
\]
is valid in a strong sense.
\item
$l$ is given explicitly by
\begin{equation} \label{lphi explicit}
l_t = \Lambda_t \circ \upsilon \text{ where }
\upsilon = (\upsilon)_{n\geq 0} \text{ and }
\upsilon_n = \tau_n \circ \varphi^{\fullstar n},
\quad n \in \ZZ_+.
\end{equation}
\end{enumerate}

\begin{rem}
In view of the injectivity of the quantum stochastic operation
(\cite{LWhom}, Proposition 2.3), Property 7$'$ implies that
\begin{equation} \label{phi to l injective}
\text{the map } \varphi \mapsto l^\varphi \text{ is injective}.
\end{equation}
\end{rem}

The next two results strengthen Property 3$'$.

\begin{propn} \label{thm: l back to k}
Let $l=l^\varphi$ and $k = k^\phi$ where
$\varphi \in CB\big(\Cil ; B(\kilhat ) \big)$
and $\phi = R_{B(\kilhat)} \varphi$. Then
$l_{t ,| \ve \ra} \in CB \big( \Cil ;| \FFock \ra \big)$ and
\[
k_{t, | \ve \ra} = R_{| \FFock \ra} l_{t, | \ve \ra} , \quad
t \in \Rplus , \ve \in \Exps .
\]
In particular, $k$ satisfies
\[
k_{t, | \ve \ra} \in CB \big( \Cil ; \Cil \ot | \FFock \ra \big),
\quad \ve \in \Exps, t \in \Rplus.
\]
\end{propn}

\begin{proof}
The first (and last) part has already been noted in Property 2$'$.
Write $\wt{k} \in \cProc (\Cil )$ for the process defined by
\[
\wt{k}_{t, | \ve \ra} =R_{| \FFock \ra} l_{t,| \ve \ra}
\in CB \big( \Cil ; \Cil \ot | \FFock \ra \big) \quad (\ve \in \Exps ,t \in \Rplus).
\]
Let $\ve = \ve (f) , \ve' = \ve (f') \in \Exps $ and $t \in \Rplus$, and consider the `form
representation' of $l$ given in Property 7$'$ and the corresponding representation of $k$.
Writing $R_\sigma$ for $R_{\Vil}$ where $\Vil = B(\kilhat)^{\ot \#
\sigma}$, Proposition~\ref{thm: R and E}
yields the identity
\[
R \big(\Omega_\sigma \circ \tau_{\# \sigma } \circ \varphi^{\fullstar \# \sigma}\big)
= \Omega_\sigma \fullcomp R_\sigma \big(\tau_{\# \sigma} \circ \varphi^{\fullstar \# \sigma}\big)
= \Omega_\sigma \fullcomp \big(\tau_{\# \sigma} \circ \phi^{\fullcomp \# \sigma}\big)
\]
$(\sigma \in \Gamma)$. Thus, integrating over $\Gamma_{[0,t]}$,
\[
R\big(l^{\ve',\ve}_t\big) = k^{\ve',\ve}_t .
\]
Therefore, using Property 5$'$,
\[
k^{\ve',\ve}_t = \cR(\lambda_{\ve'} \circ l_{t, | \ve \ra} )
=
(\id_{\Cil} \ot \lambda_{\ve'}) \circ R_{| \FFock \ra} l_{t, | \ve \ra}
= \wt{k}^{\ve', \ve}_t .
\]
The result follows.
\end{proof}

\begin{propn} \label{thm: cb l back from k}
Let $l = l^\varphi$ and $k = k^\phi$ where
$\varphi \in CB\big( \Cil ; B(\kilhat) \big)$
and $\phi = R_{B(\kilhat)} \varphi$. Then the process $l$ is completely
bounded if and only if
$k$ is, and in this case
\begin{equation} \label{k-Rl}
k_t = R_{B(\FFock)} l_t, \quad t \in \Rplus,
\end{equation}
in particular $k$ is $\Cil \ot B(\FFock)$-valued.
\end{propn}

\begin{proof}
Suppose that $l$ is completely bounded and
define the process $\wt{k}$ by $\ktilde_t = R_{B(\FFock)} l_t$.
By Proposition~\ref{thm: l back to k} and Properties 5$'$ and 5,
\[
k_{t,| \ve \ra} = R_{|\FFock \ra} (\rho_\ve \circ l_t) = (\id_{\Cil} \ot \rho_{\ve} )
\circ \ktilde_t = \ktilde_{t,| \ve \ra}
\]
$(\ve \in \Exps ,t \in \Rplus)$, and so $k$ is the completely bounded $\Cil \ot B(\FFock)$-valued
process $(R_{B(\FFock)} l_t)_{t \geq 0}$. Conversely if $k$ is completely bounded then $l$ is too,
by Property 3$'$.

\end{proof}

\section{Quantum stochastic convolution cocycles} \label{sec: E}

In this section we study quantum stochastic convolution cocycles on an
operator space coalgebra by applying the $R$-map to the theory of
quantum stochastic cocycles on an operator space (\cite{LWjfa}).
\emph{For this section an operator space coalgebra $\Cil$ is fixed.}

\begin{defn}
A completely bounded process $l\in\cProc(\Cil )$ is called a
\emph{quantum stochastic convolution cocycle}
if it satisfies
\begin{equation} \label{cbQSCC}
l_0 = \iota_{\FFock}\circ\counit \text{ and }
l_{s+t} = l_s \fullstar (\sigma_s \circ l_t) \text{ for } s,t\in\Rplus.
\end{equation}
\end{defn}
\noindent
QS convolution cocycles therefore satisfy
\begin{equation} \label{QSCC}
l_{s+t}^{\ve',\ve}
=
\la \ve'_3 , \ve_3 \ra \,
l_s^{\ve'_1, \ve_1}  \fullstar
l_t^{\ve'_2,\ve_2}
\end{equation}
for $\ve = \ve(f), \ve' = \ve (f')$ and $s,t \in \Rplus$, where
$\ve_1, \ldots ,\ve'_3$ are defined by~\eqref{epsilons}.
More generally, if $l$ is a weakly initial space
bounded process $\Cil \to \Comp$ satisfying \eqref{QSCC} then
it is called a \emph{weak} quantum stochastic convolution cocycle.
Compare this with the cocycle property for a weakly initial space
bounded process on an operator space (see Property 8 in the list
of properties of solutions of QS differential equations).

For a weak QS convolution cocycle $l$ on $\Cil$ define
\[
\lambda^{c',c}_t :=  e^{-t\la c',c\ra} \, l^{\ve',\ve}_t \text{ where }
\ve = \ve (c_{[0,t[}) \text{ and } \ve' = \ve (c'_{[0,t[})
\]
$(c, c' \in \kil, t \in \Rplus)$. Then
$\lambda^{c' ,c} := (\lambda^{c',c}_t)_{t \geq 0}$ is a convolution semigroup
and we refer to
$\{ \lambda^{c' ,c} : c, c' \in \kil \}$
as the cocycle's
\emph{associated convolution semigroups of functionals} and call $l$
\emph{Markov-regular} if $\lambda^{0,0}$ is norm continuous, in analogy to
 Markov-regular quantum stochastic cocycles
(\cite{LWjfa}).

As for standard QS cocycles, if the cocycle is contractive then
Markov-regularity implies that all of its associated convolution
semigroups of functionals are norm continuous.
Repeated application of the defining property~\eqref{QSCC} shows that, for
each $\ve = \ve(f),\, \ve' = \ve(f') \in \Exps$ and $t \in \Rplus$,
$\la\ve',\ve\ra^{-1} l^{\ve',\ve}_t$ is
the convolute of a finite number of associated convolution semigroups of
functionals of $l$. In particular
two weak QS convolution cocycles are the same if each of their
corresponding associated
convolution semigroups of functionals coincide.

\begin{lemma} \label{cocycle correspondence}
Let $l\in\cProc(\Cil)$ and $k\in\cProc(\Cil)$
be weakly initial space bounded processes related by
\begin{equation} \label{related}
k_t^{\ve' ,\ve}
 = \cR l_t^{\ve' ,\ve},
\end{equation}
for $\ve , \ve' \in \Exps, t \in \Rplus$.
Then $l$ is a weak QS convolution cocycle if and only if
$k$ is a weak QS cocycle, and in this case
$l$ is Markov-regular if and only if $k$ is.
\end{lemma}

\begin{proof}
In view of the identity
\[
\cR\big( l_s^{\ve'_1 ,\ve_1} \fullstar l_t^{\ve'_2, \ve_2} \big)
=
k_s^{\ve'_1 ,\ve_1} \circ k_t^{\ve'_2 ,\ve_2 }
\]
(in the notation~\eqref{epsilons}) the result follows from the
complete isometry of $\cR$.
\end{proof}

\begin{propn} \label{thm: l phi is a QSCC}
Let $\varphi \in CB\big(\Cil ; B(\kilhat) \big)$. Then $l^\varphi$ is a
Markov-regular weak QS convolution cocycle, each of whose convolution
semigroups of functionals is norm continuous.
\end{propn}

\begin{proof}
Let $k=k^\phi$ where $\phi = R_{B(\kilhat)} \varphi$. Then $k$ is a
Markov-regular quantum stochastic cocycle all of whose associated
semigroups are norm continuous (\cite{LWjfa}).
Since, by Proposition~\ref{thm: l back to k},
$l$ and $k$ are related by~\eqref{related} the result therefore follows from
Lemma~\ref{cocycle correspondence}.
\end{proof}

In the next section we obtain a converse by restricting to completely
positive, contractive QS convolution cocycles on a $C^*$-hyperbialgebra.
In view of the identity
\[
\int_{\Gamma_{[0,t]}} d\sigma \,
\big\la \pi_{\wh{c'}} (\sigma) ,
\varphi^{\fullstar \# \sigma} (\, \cdot \,) \pi_{\chat} (\sigma) \big\ra =
\sum_{n \geq 0}
\frac{t^n}{n!} (\omega_{\wh{c'} , \chat} \circ \varphi)^{\fullstar n} ,
\]
the convolution semigroup of functionals $\lambda^{c',c}$ associated with
the weak QS convolution cocycle $l^\varphi$ has generator
\begin{equation} \label{173}
\omega_{\wh{c'} ,\chat} \circ \varphi .
\end{equation}
This corresponds to the fact that
the semigroups associated with a Markov-regular QS cocycle $k^\phi$
on an operator space have generators
\begin{equation} \label{174}
\omega_{\wh{c'},\chat} \fullcomp \phi .
\end{equation}

Below we initiate a traffic between properties of a QS
convolution cocycle and those of its stochastic generator. Recall
Property 1$'$ for processes $l^\varphi$.
The following is easily proved either using the $R$-map, or directly.

\begin{propn} \label{realchar}
Let $l= l^\varphi$ where $\varphi \in CB \big(\Cil ; B (\kilhat) \big)$
and $\Cil$ is an operator system coalgebra.
Then
\begin{alist}
\item
$l$ is unital if and only if $\varphi (1) =0$,
\item
$l$ is real if and only if $\varphi$ is real.
\end{alist}
\end{propn}

\subsection{Opposite QS convolution cocycles}
The \emph{opposite} QS convolution cocycle relation, for processes in
$\cProccb (\Cil )$, is
\[
l_0 = \iota_{\FFock}\circ\counit \text{ and }
l_{s+t} = (\sigma_s \circ l_t)  \fullstar l_s,
\]
which involves the natural identifications
$B(\FFock_{[s,s+t[}) \olot B(\FFock_{[0,s[}) = B(\FFock_{[0,s+t[})$,
for $s,t\in\Rplus$. Completely bounded processes which satisfy a QS
differential equation of the form
\[
dl_t = d\Lambda_{\varphi}(t) \star l_t, \;\; l_0 = \iota_{\FFock} \circ \Cou,
\]
for $\varphi \in CB(\Cil;B(\kilhat))$, are opposite QS convolution
cocycles; they are given explicitly by
\[
l_t = \Lambda_t \circ \upsilon \text{ where }
\upsilon_n = \varphi^{\star n}, n\in\ZZ_+
\]
(cf.\ Properties 6$'$ and 7$'$ in Section~\ref{sec: D}), with
${}^\varphi l$ being an appropriate notation.
There is a bijective correspondence between the set of QS convolution
cocycles treated in this paper and the set of opposite QS convolution
cocycles. This is effected by time-reversal, as in~\cite{LWjfa}. Opposite
QS convolution cocycles have convolution semigroup representation as QS
convolution cocycles do, but with the semigroups appearing in the reverse
order. In particular time-reversal exchanges $l^\varphi$ and
${}^\varphi l$. One may also view the correspondence in terms of the
opposite coproduct $\Com^{\opposite}$.

In~\cite{LSqscc1} we actually worked with opposite cocycles (thus the
convolvands in (5.1) and the $a_{(i)}$'s in (5.2), on p.\ 595 of that
paper, should both have appeared in the reverse order, with the notation
${}^\varphi l$ being more appropriate for the opposite QS convolution
cocycles generated there). The results of that paper are equally valid
for QS convolution cocycles on coalgebras defined as here through the
relations~\eqref{cbQSCC} and~\eqref{QSCC}.

\section{Completely positive QS convolution cocycles} \label{sec: F}

In this section we characterise the Markov-regular QS convolution cocycles amongst the completely positive and
contractive processes on on
a $C^*$-hyperbialgebra $\El$, as those which satisfy a coalgebraic
quantum stochastic differential equation with completely bounded
coefficient of a particular form. We also give the general form of the coefficient of the QS
differential equation. Recall Theorem~\ref{CPstandard} and
the notations~\eqref{Delta QS}.

\begin{thm} \label{genstruct}
Let $\elg$ be a $C^*$-hyperbialgebra and let $l\in \cProc(\elg)$.
Then the following are equivalent\tu{:}
\begin{rlist}
\item
$l$ is a Markov-regular, completely positive, contractive QS convolution
 cocycle\tu{;}
\item
$l=l^{\varphi}$ where $\varphi \in CB \big( \elg ; B(\kilhat) \big)$
satisfies $\varphi (1) \leq 0$
and may be decomposed as follows\tu{:}
\begin{equation}    \label{CPC1}
   \varphi  =  \psi  -  \counit ( \cdot )
   \left(\Projk + | e_0 \ra \la \chi | + | \chi \ra \la e_0 | \right)
   \end{equation}
for some completely positive map $\psi: \elg \to B(\kilhat)$ and vector
$\chi \in \kilhat$\tu{;}
\item
there is a *-representation $(\rho,\Kil)$ of $\El$, a
contraction $D\in B(\kil;\Kil)$ and a vector $\xi\in\Kil$,
such that $l=l^\varphi$ where
\begin{equation} \label{alternative CPC}
\varphi (x) =
\begin{bmatrix}\la\xi | \\ D^*\end{bmatrix}
\big(\rho(x) - \counit (x) I_{\Kil}\big)
\begin{bmatrix}|\xi \ra & D\end{bmatrix}
+ \counit (x)\varphi (1)
\end{equation}
$(x\in\El)$, and $\varphi (1)$ is nonpositive with block matrix of the
form
\[
\begin{bmatrix}
* & * \\
* & D^*D - I_{\kil}
\end{bmatrix}.
\]
\end{rlist}
\end{thm}

\begin{proof}
For the proof of the equivalence of (i) and (ii)
we may suppose that $\El$ is faithfully and nondegenerately represented in
$B(\init)$, say, in such a way that the counit extends to a normal state
$\counit''$ on $\El''$.
(This may be achieved
by taking the direct sum of an arbitrary faithful nondegenerate
representation and the GNS representation
$(\init_{\counit} , \pi_{\counit} , \xi_{\counit})$, so that $\counit$
is extended by the vector state $\omega_{(0, \xi_{\counit})}$.
Alternatively, take the universal representation
and bidual map $\counit^{**}$.) Note that $\counit''$ is necessarily
*-homomorphic.

Suppose first that (i) holds and let $\{ \gamma_{c',c} :c',c \in \kil\}$
be the generators
of the associated convolution semigroups of functionals of $l$.
Let $k \in \Procdagger (\Eil \Pto \Eil)$
be the process $(R_{B(\FFock)} l_t)_{t \geq 0}$. By
Proposition~\ref{Requiv},
$k$ is completely positive and contractive. Moreover~\eqref{related} holds
so that $k$ is a
Markov-regular QS cocycle on $\Eil$. In view of Theorem~\ref{CPstandard},
it follows that $k=k^\phi$ for
some map $\phi$ of the form~\eqref{CPC0} with
constituents $\Psi \in CP\big(\Eil ; \Eil''\, \olot\, B(\kilhat)\big)$ and
$J \in \Eil''\, \olot\, | \kilhat \ra$,
say.   Thus, from the definition of $k$,
\begin{equation} \label{phi and gamma}
\omega_{\wh{c'} , \chat} \fullcomp \phi = ( \id_{\El} \ot \gamma_{c',c} ) \circ \Delta
\end{equation}
$(c,c' \in \kil)$.
Now define
$\varphi , \psi \in CB \big( \El; B(\kilhat)\big)$ and $\chi \in \kilhat$ by
\[
\varphi = \counit \fullcomp \phi , \quad
\psi = ( \counit''\, \olot\, \id_{B(\kilhat)} ) \circ \Psi  \text{
and }
\la \chi | = (\counit''\, \olot\, \id_{\la \kilhat |} ) (J),
\]
noting that, by the complete positivity of $\counit$,
$\varphi (1) \leq 0$ and $\psi$ is completely positive.
We claim that $l = l^\varphi$ and that $\varphi$ has the
decomposition~\eqref{CPC1}.
By~\eqref{phi and gamma},
$\omega_{\wh{c'} , \chat} \circ \varphi = \gamma_{c',c}$ and so,
by~\eqref{173}, the QS convolution cocycles $l^\varphi$ and $l$ have the
same associated convolution semigroups and are therefore equal.
In view of the multiplicativity of $\counit''$,
\[
(\counit'' \olot \id_{B(\kilhat)})\big( (x \ot | e_0 \ra ) J\big)
= \counit (x) | e_0 \ra \la \chi |,
\]
for $x\in\El$.
Now, using the fact that $\counit''$ is real to obtain the adjoint
identity, collecting terms yields the decomposition~\eqref{CPC1}, so (ii)
holds.

Suppose conversely that (ii) holds. As in the proof of
Proposition~\ref{thm: l phi is a QSCC}, let $k=k^\phi$ where
$\phi =
R_{B(\kilhat)} \varphi \in CB \big( \Eil ; \Eil \ot B(\kilhat )\big)$.
Then $\phi$ has the form~\eqref{CPC0}, with $\Psi = R_{B(\kilhat)} \psi$ and
$J = I_{\init} \ot \la \chi |$, moreover it follows from
Proposition~\ref{Requiv} that $\Psi$ is
completely positive and $\phi (1) \leq 0$. Thus, by
Theorem~\ref{CPstandard}, the Markov-regular weak QS cocycle $k$ is
completely positive and contractive. Therefore, by
Proposition~\ref{thm: cb l back from k},
Proposition~\ref{Requiv} and
Lemma~\ref{cocycle correspondence}, (i) holds.

Again suppose that (ii) holds. Let
\begin{equation} \label{psi Stinespring}
\begin{bmatrix} \la\xi | \\ D^* \end{bmatrix}
\rho(\cdot)
\begin{bmatrix} |\xi\ra & D \end{bmatrix}
\end{equation}
be a minimal Stinespring decomposition of $\psi$. Thus $(\rho,\Kil)$ is a
unital $C^*$-representation of $\El$, $\xi$ is a vector in $\Kil$, $D$
is an operator in $B(\kil;\Kil)$ (and~\eqref{minimal} below holds).
Identity~\eqref{alternative CPC} follows, with
\[
\varphi (1) =
\begin{bmatrix}
\|\xi\|^2 - 2\re \alpha & \la D^*\xi - c| \\
| D^*\xi - c\ra & D^*D - I_{\kil}
\end{bmatrix},
\]
where $\binom{\alpha}{c}=\chi$, so (iii) holds.

Conversely, suppose that (iii) holds. Then, writing
\[
\begin{bmatrix}
t & \la d | \\
|d\ra & D^*D - I_{\kil}
\end{bmatrix}
\]
for the block matrix form of $\varphi (1)$, $\varphi$ has the
form~\eqref{CPC1} where $\psi$ is given by~\eqref{psi Stinespring} and
\[
\chi = \binom{\frac{1}{2}(\|\xi\|^2 - t)}{D^*\xi - d}
\]
so (ii) holds. This completes the proof.
\end{proof}

\begin{rems}
An alternative proof of the above theorem, which directly establishes the
equivalence of (i) and (iii) without appeal to Theorem~\ref{CPstandard}
on standard QS cocycles (whose proof depends on the Christensen-Evans
Theorem), is given in~\cite{Sthesis}.

In (iii) the following
\emph{minimality} condition on the quadruple $(\rho,\Kil,D,\xi)$ may be
assumed:
\begin{equation} \label{minimal}
\rho(\El)\big(\Comp\xi + \Ran D\big) \text{ is dense in } \Kil.
\end{equation}
Under minimality there is uniqueness too: if $(\rho',\Kil',D',\xi')$ is
another quadruple as in (iii) then there is a unique isometry
$V\in B(\Kil;\Kil')$ (unitary if this quadruple is also minimal)
satisfying
\[
VD=D', V\xi =\xi' \text{ and } V\rho (x)=\rho'(x)V \text{ for } x\in\El.
\]

By a characterisation of nonnegative block matrix operators
(see, for example, Lemma 2.2 in~\cite{GLSW}), if $\varphi$
is the stochastic generator of a Markov-regular, completely positive,
contractive QS convolution cocycle then $\varphi (1)$ has the form
\[
\begin{bmatrix}
t & \la C^{1/2}e| \\
|C^{1/2}e\ra & -C
\end{bmatrix}
\]
for a nonnegative contraction $C$, a unique
vector $e\in\ol{\Ran}\, C$ and a real number $t$
satisfying $t\leq -\| e\|^2$. Moreover, with
respect to any decomposition~\eqref{alternative CPC}, $C=I_{\kil}-D^*D$.

Unitality for the cocycle is equivalent to its stochastic generator being
expressible in the form
\[
\begin{bmatrix} \la\xi | \\ D^* \end{bmatrix}
\big(\rho - \iota_{\Kil}\circ\counit\big) (\cdot )
\begin{bmatrix} |\xi\ra & D \end{bmatrix}
\]
where $D$ is isometric and $\rho (1)$ coincides with the identity operator
on $\Comp\xi + \Ran D$.

It also follows from the above proof that if $k$ is the QS cocycle on a
$C^*$-hyperbialgebra $\El$,
given by $~k_t = R_{B(\FFock)} l_t$ where $l$ is a Markov-regular,
completely positive,
contractive QS convolution cocycle on $\El$, then $k = k^\phi $ where
the stochastic generator $\phi$ is expressible in the form
\[
x \mapsto \psi (x) - x \ot \big( \Delta^{QS} + | \chi \ra \la e_0 |
+ | e_0 \ra \la \chi | \big)
\]
for some completely positive map $\psi : \El \to \El \ot B(\kilhat)$ and
vector $\chi \in \kilhat$.
Note that no appeal to a concrete realisation of $\El$ is needed in this
decomposition.
\end{rems}

\section{Homomorphic QS convolution cocycles} \label{sec: G}

In this section we characterise the stochastic generators of
Markov-regular *-homomorphic convolution cocycles on a $C^*$-bialgebra,
by applying the $R$-map to the characterisation of the generators of
Markov-regular multiplicative cocycles obtained in~\cite{LWhom}.
Thus let $\Bil$ be a $C^*$-bialgebra.

\emph{Weak multiplicativity} for a process
$l\in\cProcdagger (\Bil )$ is the following property:
\[
l^{\ve' ,\ve}_t (x^*y) = l^\dagger_{t, |\ve'\ra}(x)^* l_{t,|\ve\ra}(y)
\]
$(\ve , \ve' \in \Exps , x,y\in\Bil, t \in \Rplus )$. If the
$C^*$-bialgebra is concretely realised on a Hilbert space then weak
multiplicativity for a process
$k \in \Procdagger (\Bil \Pto \Bil)$ reads
\[
k_t (x^*y) = k^\dagger_t (x)^* k_t(y)
\]
$(x,y \in \Bil, t \in \Rplus )$, an identity in Hilbert space operators.
In view of the remark at the end of the introduction, if
$k \in \Procdagger (\Bil \Pto \Bil)$ is both weakly multiplicative and
real then it is bounded,
and so *-homomorphic --- in particular it is completely bounded.

\begin{thm}
Let $l=l^\varphi$ where $\varphi \in CB\big(\Bil ;B (\kilhat) \big)$.
Then the following are equivalent
\begin{alist}
\item
$l$ is weakly multiplicative\tu{;}
\item
$\varphi$ satisfies
\begin{equation} \label{varphi structure}
\varphi (xy) =
\varphi (x) \counit (y) +
\counit (x) \varphi(y) +
\varphi(x) \Projk \varphi (y)
\end{equation}
\tu{(}$x,y\in\Bil$\tu{)}.
\end{alist}
\end{thm}

\begin{proof}
For the proof we may suppose without loss of generality that the
$C^*$-bialgebra $\Bil$ is concretely realised, in $B(\init)$ say. Let
$\phi = R_{B(\kilhat)} \varphi$ and set $k = k^\phi$. Since
$\Ran \phi \subset \Bil \ot B(\kilhat)$,
Theorem 3.4 and Corollary 4.2 of~\cite{LWhom} imply that $k$ is weakly multiplicative if and
only if $\phi$ satisfies
\begin{equation} \label{phi structure}
\phi (xy) = \phi (x) \iota_{\kilhat}(y) +
\iota_{\kilhat}(x) \phi (y) + \phi (x)
\big(I_{\init}\ot\Projk\big) \phi (y).
\end{equation}

If~\eqref{phi structure} holds then, applying the homomorphism $\counit
\ot \id_{B(\kilhat)}$
to both sides yields~\eqref{varphi structure}. Conversely, suppose
that~\eqref{varphi structure}
holds and set $\wt{\varphi} = \id_{\Bil} \ot \varphi$ and $\wt{\counit} =
\id_{\Bil} \ot (\iota_{\kilhat}\circ\counit )$,
so that
\begin{equation} \label{phi tilde delta}
\phi = \wt{\varphi} \circ \Com \text{ and }
\wt{\counit} \circ \Com = \iota_{\kilhat}.
\end{equation}
Then, for simple tensors $X = x_1 \ot x_2$ and $Y=y_1 \ot y_2$ in $\Bil
\ot \Bil$,
\[
x_1y_1 \ot \varphi(x_2y_2) = x_1y_1 \ot \big\{ \varphi (x_2) \counit (y_2)
+ \counit(x_2) \varphi(y_2)
    + \varphi(x_2)\Projk \varphi (y_2) \big\} ,
\]
or
\[
\wt{\varphi} (XY) =
\wt{\varphi} (X) \wt{\counit} (Y) + \wt{\counit}(X) \wt{\varphi}(Y)
  + \wt{\varphi} (X) \big(I_{\init}\ot\Projk\big)  \wt{\varphi} (Y) .
\]
By linearity and continuity this holds for all $X,Y \in \Bil \ot \Bil$, in
particular, for $X = \Com x$ and $Y= \Com y$. Therefore,
by~\eqref{phi tilde delta} and the multiplicativity of $\Com$,
\eqref{phi structure} holds.

It therefore remains only to show that $l$ is weakly multiplicative if and
only if $k$ is. Recall that $l\in\cProcdagger (\Bil)$ and
$l^\dagger = l^\psi$ where $\psi = \varphi^\dagger$.
Let $u,u' \in \init$, $t \in \Rplus$ and $\ve, \ve' \in \Exps$. If
$l$ is
weakly multiplicative then
$\big\la \ve', l_t (xy) \ve \big\ra =
\big\la l^\dagger_t (x^*) \ve', l_t (y) \ve \big\ra$, so
\begin{equation}
\label{l(xy) to k(xy)}
\big\la u'\ve', (\id_{\Bil} \odot l_t )(XY) u \ve \big\ra =
\big\la (\id_{\Bil} \odot l^\dagger_t) (X^*) u'\ve',
(\id_{\Bil} \odot l_t) (Y) u\ve \big\ra
\end{equation}
holds for all $X,Y \in \Bil \odot \Bil$. Now it follows, from the identity
\[
(\id_{\Bil} \odot l^\#_t) (X) u \ve =
\big(\id_{\Bil} \ot l^\#_{t,| \ve\ra}\big) (X) u
\]
(where $l^\#$ stands for $l$ or $l^\dagger$) and Property 2$'$,
that both sides of~\eqref{l(xy) to k(xy)} are continuous in both $X$
and $Y$, giving an identity for all $X,Y \in\Bil \ot \Bil$. Setting
$X= \Com x$ and $Y= \Com y$ and using the multiplicativity of $\Com$, this
identity becomes a statement of the weak multiplicativity of $k$.

Suppose conversely that $k$ is weakly multiplicative. Set $k^\dagger =
k^\psi$ for
$\psi = \phi^\dagger$, and let $x,y \in \Bil$. First note that the
identity
\[
(\counit \ot \id_{\la \FFock|}) (X^*) (\counit \ot \id_{| \FFock \ra} (Y)
 = \counit (X^*Y)
\]
is obvious for $X,Y \in \Bil \odot | \FFock \ra$ and so holds for
$X,Y \in \Bil \ot | \FFock \ra$ by continuity. Set
$X=k^\dagger_{t,| \ve \ra} (x^*)$ and $Y=k_{t,| \ve' \ra} (y)$. Then
$X,Y \in \Bil \ot | \FFock \ra$ and so
\begin{align*}
\big\la l^\dagger_t (x^*) \ve, l_t(y)\ve' \big\ra
  &= (\counit \ot \id_{|\FFock \ra})(X)^* (\counit \ot \id_{| \FFock \ra}) (Y) \\
  &= \counit (X^*Y) \\
  &= (\counit \circ \omega_{\ve,\ve'} \fullcomp k_t) (xy) \\
  &= (\omega_{\ve,\ve'} \circ \counit \fullcomp k_t) (xy) = \big\la \ve
,l_t(xy) \ve' \big\ra .
\end{align*}
Thus $l$ is weakly multiplicative. This completes the proof.
\end{proof}

Combining this result with Proposition~\ref{realchar},
Theorem~\ref{genstruct} and Theorem~\ref{A6} we obtain the advertised
characterisation of the stochastic generators of Markov-regular
*-homomorphic convolution cocycles on a $C^*$-bialgebra.

\begin{thm}\label{*homchar}
Let $\Bil$ be a $C^*$-bialgebra and let $l \in \cProc (\Bil )$.
Then the following are equivalent\tu{:}
\begin{rlist}
\item
$l$ is a Markov-regular, *-homomorphic \tu{(}and unital\tu{)} QS
convolution cocycle on $\Bil$\tu{;}
\item
$l=l^\varphi$ where $\varphi \in CB \big( \Bil ; B(\kilhat) \big)$
satisfies
\begin{equation}  \label{hom structure}
\varphi(x^*y) =
\varphi (x)^* \counit (y) +
\counit (x)^* \varphi(y) +
\varphi (x)^* \Projk\varphi(y) \;\; \big(\text{and } \, \varphi(1)=0\big) ;
\end{equation}
\item
there is a vector $c \in \kil$ and \tu{(}unital\tu{)} *-homomorphism
$\pi : A \to B(\kil)$ such that $l = l^\varphi$ where
\begin{equation} \label{inner generator}
\varphi (x) = \begin{bmatrix} \la c | \\ I_{\kil} \end{bmatrix}
\big( \pi (x) - \counit (x) I_{\kil} \big)
\begin{bmatrix}|c \ra & I_{\kil} \end{bmatrix},
\quad x \in \Bil .
\end{equation}
\end{rlist}
\end{thm}

\begin{rem}
In fact, as is shown in the appendix, the relation~\eqref{hom structure}
for a linear map $\varphi$ (an \emph{$\counit$-structure map} in the
terminology used there) entails the implemented form~\eqref{inner
generator}, in particular the complete boundedness of $\varphi$.
\end{rem}

The characterisations of stochastic generators of completely positive,
contractive QS convolution cocycles and *-homomorphic QS convolution
cocycles in Theorems~\ref{genstruct} and~\ref{*homchar} may be used to
derive dilation theorems for QS convolution cocycles (see~\cite{Sthesis}),
of the type obtained for standard QS cocycles in~\cite{GLW}
and~\cite{GLSW}. These characterisations are also used to establish the
main result in~\cite{FrS}, that every Markov-regular Fock space quantum
L\'evy process can be realised as a limit of (suitably scaled) random
walks.

\section{Axiomatisation of topological quantum L\'evy processes} \label{Levy}

Defining quantum L\'evy process on a $C^*$-bialgebra requires certain
modifications of the original, purely algebraic, definition of Accardi,
Sch\"urmann and von Waldenfels (\cite{asw}, \cite{Schurmann}).
The problem is how to build convolution increments of the process given
that, in general, multiplication $\blg \odot \blg \to \blg$ need not
extend continuously to $\blg \ot \blg$. (This is a commonly met difficulty
in the theory of topological quantum groups, see~\cite{kus}). Below we
outline two ways of overcoming this obstacle.

The simplest idea is to define a quantum L\'evy process using only the
concept of distributions.

\begin{defn}  \label{wLp}
A \emph{weak quantum L\'{e}vy process} on a $C^*$-bialgebra $\blg$
over a unital *-algebra-with-state
$(\Alg , \omega)$ is a family
$\big(j_{s,t}\! :\blg \Pto \Alg\big)_{0 \leq s \leq t}$
of unital *-homomorphisms such that
the functional  $\lambda_{s,t}:= \omega \circ j_{s,t}$ is continuous and
satisfies the following conditions,
for $0\leq r \leq s \leq t$:

\begin{enumerate}
[{(wQLP}i)]
\item       \label{wQLPi}
$\lambda_{r,t} = \lambda_{r,s} \star \lambda_{s,t}$\tu{;}
\item  \label{wQLPii}
$\lambda_{t,t} = \Cou $\tu{;}
\item \label{wQLPiii}
$\lambda_{s,t} = \lambda_{0,t-s}$\tu{;}
\item \label{wQLPiv}
\[
\omega \left( \prod^n_{i=1} j_{s_i,t_i} (x_i) \right) = \prod^n_{i=1}
\lambda_{s_i,t_i} (x_i)
\]
whenever $n \in \bn$, $x_1, \ldots, x_n \in\blg$ and
the intervals $[s_1,t_1[,\ldots ,[s_n, t_n[$ are disjoint\tu{;}
\item  \label{wQLPv}
$\lambda_{0,t}  \to \Cou$ pointwise as $t \to 0$.
\end{enumerate}
A weak quantum L\'{e}vy process on a $C^*$-bialgebra $\blg$ is called
\emph{Markov-regular} if
$\lambda_{0,t} \to\Cou$ in norm, as $t \to 0$.
\end{defn}

The family $\lambda:=\big(\lambda_{0,t}\big)_{t \geq 0}$ is a pointwise
continuous convolution semigroup of functionals on $\blg$, called the
\emph{one-dimensional distribution} of the
process; if the process is Markov-regular then $\lambda$ has a convolution
generator which is also referred to as the \emph{generator} of the weak
quantum L\'{e}vy process.
Two weak quantum L\'{e}vy processes on $\blg$,
$j^1$ over $(\Alg^1 ,\omega^1)$ and $j^2$ over $(\Alg^2 ,\omega^2)$,
are said to be \emph{equivalent} if they satisfy
\[
\omega^1 \circ j^1_{s,t}  = \omega^2 \circ j^2_{s,t}
\]
for all $0\leq s\leq t$, in other words if their one-dimensional
distributions coincide; if they are Markov-regular then this is
equivalent to equality of their generators.

\begin{rems}
Note that the above definition of a weak quantum L\'evy process, in
contrast to the definition
of a quantum L\'evy process on an algebraic *-bialgebra, does not yield a
recipe for expressing the joint moments of the process increments
corresponding to overlapping time intervals, such as
\[
\omega(j_{r,t}(x) j_{s,t} (y)) \text{ where } r,s < t.
\]
To achieve the latter, one would have to formulate the weak convolution
increment property (wQLP\ref{wQLPi}) in greater generality and assume
certain commutation relations between the increments corresponding to
disjoint time intervals.
For other investigations of the notion of independence in noncommutative
probability, in the absence of commutation relations being imposed, we
refer to the recent paper \cite{hkk}.
\end{rems}

As in the algebraic case, the generator of a Markov-regular weak quantum
L\'evy process vanishes on $1_{\blg}$, is real and is \emph{conditionally
positive}, that is positive on the kernel of the counit.
Observe that if $l\in \cProc(\blg)$ is a unital *-homomorphic QS
convolution cocycle then, defining
$\Alg := B(\FFock)$, $\omega:= \omega_{\ve(0)}$, and
$j_{s,t} := \sigma_s \circ l_{t-s}$ for all $0\leq s\leq t$,
we obtain a weak quantum L\'evy process on $\blg$, called a \emph{Fock
space quantum L\'evy process}, Markov-regular if $l$ is.
The proof of the following theorem closely mirrors the proof of
Sch\"urmann's reconstruction theorem for the purely algebraic case
(\cite{Schurmann}, see also \cite{LSqscc1}); all the necessary continuity
properties follow from the results in the appendix.

\begin{thm} \label{Crecon}
Let $\gamma$ be a real,
conditionally
positive linear functional on $\blg$ vanishing at $1_{\blg}$.
Then there is a \tu{(}Markov-regular\tu{)}
Fock space quantum L\'{e}vy process with generator $\gamma$.
\end{thm}

\begin{proof}
The proof uses a GNS-style construction. Let $D= \Ker \Cou \big/ N$ where
$N$ is the following subspace of $\Ker \counit$:
\[
\big\{ x \in \Ker \epsilon \, \big| \; \gamma (x^* x) =0 \big\}.
\]
Then $\big([x],[y]\big) \mapsto \gamma (x^*y)$ defines an inner product
on $D$. Let $\kil$ be the Hilbert space completion
of $D$. The prescription $\pi (x): [z] \mapsto [xz]$ defines bounded
operators on $D$, whose extensions make
up a unital representation of $\blg$ on $\kil$ satisfying
\[
\big\la \pi (x)[y],[z] \big\ra = \big\la [y], \pi(x^*) [z] \big\ra .
\]
Furthermore the linear map
$\delta : x \mapsto  |d(x)\ra $, where $d(x)=[x-\counit (x) I_{\kil}]$,
is easily seen to be a
($\pi$-$\epsilon$)-derivation $\blg\rightarrow |\kil\ra $ satisfying
\[
\delta (x)^* \delta (y) =
\gamma (x^*y) - \gamma (x)^*\epsilon (y) - \epsilon (x)^* \gamma (y).
\]
Theorem \ref{A6} therefore implies that the map
$\varphi:\blg \rightarrow B(\kilhat)$, with block matrix form given by the
prescription \eqref{the new 6.6} (with
$\lambda = \gamma$ and $\chi = \counit$), is completely bounded.
Setting $l = l^\varphi$, Theorem~\ref{*homchar} implies that the
Markov-regular weak QS convolution cocycle $l$ is unital and
*-homomorphic. Since $\varphi^0_0 =\gamma$ the result follows.
\end{proof}

\begin{cor}
Every Markov-regular weak quantum L\'{e}vy process is equivalent to a Fock
space quantum L\'{e}vy process.
\end{cor}

Another notion, in a sense intermediate between weak quantum L\'{e}vy
processes and Fock space quantum L\'{e}vy processes,
can be formulated in terms of product systems --- a similar idea
is mentioned in a recent paper of Skeide (\cite{Skeiden}). Recall that a
\emph{product system of Hilbert spaces} is a `measurable' family of
Hilbert spaces $E=\{E_t: t \geq 0 \}$, together with unitaries
$U_{s,t}:E_s \ot E_t \to E_{s+t}$ ($s,t \geq 0$) satisfying associativity
relations:
\begin{equation} \label{ass}
U_{r+s, t} (U_{r,s} \ot I_{t}) = U_{r, s+t} (I_{r} \ot U_{s,t})
\end{equation}
($r,s,t \in \br_+$), where $I_s$ denotes the identity operator on $E_s$.
A \emph{unit} for the product system $E$ is a `measurable' family
$\{u(t): t \geq0\}$ of vectors with $u(t) \in E_t$ and
$u(s+t) = U_{s,t} \big(u(s) \ot u(t)\big)$ for all
$s,t \geq 0$ (the unit is \emph{normalised} if, for all $t \geq 0$,
$\|u(t)\|=1$). For the precise definition we refer to \cite{Arveson}.
The unitaries $U_{s,t}$ implement isomorphisms
$\sigma_{s,t}: B(E_s \ot E_t) \to B(E_{s+t})$.

\begin{defn}
A \emph{product system quantum L\'evy process} on $\blg$
over a product-system-with-normalised-unit $(E,u)$ is a family
$\big(j_t: \blg \rightarrow B(E_t)\big)_{ t \geq 0}$ of unital
*-homomorphisms satisfying the following conditions:
\begin{rlist}
[(psQLPi)] \label{EuQLPi}
\item
$j_{r+s} = \sigma_{r,s}\circ \big(j_r\fullstar j_s\big)$,
\item [(psQLPii)] \label{EuQLPii}
$j_0 = \iota_0\circ\Cou$,
\item [(psQLPiii)] \label{EuQLPiii}
$\omega_{u(t)}\circ j_t \to\Cou$
pointwise as $t \to 0$,
\end{rlist}
for $r,s \geq 0$,
where $\iota_0$ denotes the ampliation $\Comp\to B(E_0)$.
\end{defn}

The `exponential' product system is given by
$E_t = \FFocktot$ and $U_{s,t}= I_s\ot S_{s,t}$ where $S_{s,t}$ denotes
the natural shift $\FFocktot\to\FFock_{[s,s+t]}$ and the exponential
property of symmetric Fock space is invoked. Clearly every Fock space
quantum L\'{e}vy process may be viewed as a  product system quantum
L\'{e}vy process over $(E,\Omega)$ where $\Omega$ is the normalised unit
given by $\Omega (t)= \ve(0)\in\FFock_{[0,t[}$, $t\geq 0$.

\begin{propn}
Each product system quantum L\'{e}vy process on $\blg$ naturally
determines a weak quantum L\'{e}vy process on $\blg$ with the same
one-dimensional distribution.
\end{propn}

\begin{proof}
Let $j$ be a quantum L\'evy process on $\Bil$ over a
product-system-with-normalised-unit $(E,u)$.
We use an inductive limit construction.
Define $\wt{\Alg} := \bigcup_{t \geq 0} (B(E_t),t)$ and introduce on
$\wt{\Alg}$ the relation: $(T,r) \equiv (S,s)$ if
there is $t \geq \max\{r,s\}$ such that
$\sigma_{r,t-r} (T \ot I_{t-r} )= \sigma_{s,t-s} (S \ot I_{t-s})$,
in other words we identify operators with common ampliations.
The associativity relations~\eqref{ass} imply that
$\equiv$ is an equivalence relation. Define $\Alg = \wt{\Alg}/\!\equiv$
and introduce the structure of a unital *-algebra on $\Alg$,
consistent with the pointwise operations:
\[
(T,t) + (S,t) = (T+S,t), \;\;
(S,t) \cdot (T,t) = (ST, t), \;\;
(T,t)^* = (T^*,t)
\]
($t \geq 0, S,T \in B(E_t)$).
The map $\wt{\omega}:\wt{\Alg} \to \bc$ defined by $\wt{\omega} (T,t) = \omega_{u(t)}(T)$
induces a state $\omega$ on $\Alg$.
For $s,t\in\Rplus$ define
\[
j_{s,t}:\Bil\to\Alg \text{ by }
x\mapsto \left[\sigma_{s, t-s} (I_s \ot j_{t-s}(x) ) \right]_{\equiv}.
\]
It is easy to see that the family
$\big(j_{s,t}\big)_{0 \leq s \leq t}$ is a weak quantum
L\'evy process on $\blg$ over $(\Alg, \omega)$.
\end{proof}

The construction in the above proof, informed by the case of QS
convolution cocycles, is a special case of the familiar construction of
$C^*$-algebraic inductive limits. The completion of $\Alg$ with respect to
the norm induced from $\wt{\Alg}$ is a unital $C^*$-algebra that may be
called the \emph{ $C^*$-algebra of finite range operators on the product
system $E$}.

\begin{rem}
A form of reconstruction theorem also holds for completely positive QS convolution cocycles.  It is easily seen that if
$l\in\cProc(\elg)$ is a Markov-regular, unital, completely
positive QS convolution cocycle on a $C^*$-hyperbialgebra $\elg$, then the generator of its
Markov convolution semigroup is real, vanishes at $1_{\elg}$ and is conditionally positive. The GNS-type
construction from the proof of Theorem \ref{Crecon} yields a completely bounded map $\varphi:\alg \to B(\kilhat)$ for
which the cocycle $l^{\varphi}$ is unital and completely positive according to Proposition \ref{realchar} and Theorem \ref{genstruct}
(of course there is no reason why it should be *-homomorphic, if $\elg$ is not a
$C^*$-bialgebra).
Clearly the Markov convolution semigroup of $l^{\varphi}$ coincides with that of $l$.
\end{rem}

\section{Examples} \label{section: H}

In this section we consider *-homomorphic convolution
cocycles on three types of $C^*$-bialgebra, namely algebras of continuous functions on compact semigroups, universal $C^*$-algebras of discrete
groups, and full compact quantum groups. We focus on connections between the results obtained in this paper and
the case of purely algebraic convolution cocycles analysed in its predecessor, \cite{LSqscc1}. Recall that
in \cite{LSqscc1} the basic object is an algebraic *-bialgebra (or even
coalgebra) $\Blg$, and coalgebraic
QS differential equations are driven by coefficients in $L(\Blg; \Opdagger(\Dhat))$, where $D$ is some dense
subspace of the noise dimension space $\kil$. Processes $\Vlg\to\Comp$, now for a vector space $\Vlg$, are families
$k=(k_t)_{t \geq 0}$ of maps
$\Vlg \to \Op(\ExpsD)$;
we denote the space of these by $\cProc(\Vlg :\ExpsD)$,
and write $\cProcdagger(\Blg :\ExpsD)$ for the subspace of
$\Opdagger(\ExpsD)$-valued processes.
\emph{Pointwise H\"older-continuity} for such a process $k$ means that
each of the vector-valued functions $t \mapsto k_t(x) \ve$ should be locally
H\"{o}lder-continuous with exponent  $1/2$. Note that it is a
weaker form of continuity than the one that arises when $\Vlg$ is an
operator space (cf.\ Properties
\ref{stHol1} and \ref{stHol2}$'$
after Theorems \ref{thm: existence} and \ref{thm: coalg existence}).

The notation introduced after Theorem \ref{thm: coalg existence} extends
as follows: for $\varphi \in L(\Blg; \Op(\Dhat))$,
\eqref{lphi explicit} still
defines a process $\cProc(\Blg :\Exps_D)$
(again written $l^{\varphi}$) which (uniquely) satisfies the QS
differential equation~\eqref{OScoalgQSDE}, now understood in the sense
of~\cite{LSqscc1}, and is a QS convolution cocycle with respect to the
purely algebraic coalgebra structure. If the coefficient $\varphi$
lies in $L(\Blg; \Opdagger(\Dhat))$ then
$l^{\varphi}\in\cProcdagger(\Blg :\Exps_D)$.

\subsection{Commutative case: continuous functions on a semigroup }

Let $H$ be a compact semigroup with identity $e$ and let
$\Bil$ denote $C(H)$, the algebra of continuous complex-valued functions
on $H$.
Then $\Bil$ has the structure of a $C^*$-bialgebra with comultiplication
and counit given
by
\[
\Com(F) (h,h') = F(hh') \;\;\
\text {and }
\Cou (F) = F(e)
\]
($h,h' \in H, F\in \Bil$),
courtesy of the natural identification $\Bil\ot \Bil \cong C(H\times H)$.

Following standard practice in quantum probability
(going back to \cite{AFL} and beyond), any $H$-valued stochastic
process $X=\big(X_t\big)_{t\geq 0}$
on the probability space $(\Omega,\mathfrak{F},\Prob)$, may be described by
a family of unital *-homomorphisms $\big(l_t\big)_{t \geq 0}$ given by
\[
l_t :\Bil\to L^{\infty}(\Omega, \mathfrak{F}, \Prob),
\quad
F \mapsto F \circ X_t,
\]
in turn these homomorphisms uniquely determine the original process.

Recall that a process $X$ on a semigroup with identity is called a
\emph{L\'evy process} if
it has identically distributed, independent increments,
$\Proc(\{X_0=e\})=1$ and
the distribution of $X_t$ converges weakly to the Dirac measure $\delta_{\{e\}}$
(the distribution of $X_0$) as $t$ tends to $0$.
In general every L\'evy process on a semigroup may be equivalently realised, in the sense
of equal finite-dimensional distributions
(see \cite{Schurmann}, \cite{LSqscc1}), as a quantum L\'evy process on a *-bialgebra
(\cite{Schurmann}, \cite{FranzSchott}).

As is well known, not all L\'evy processes have stochastic
generators defined on the whole of $\Bil$. In our language, this
corresponds to the fact that not all *-homomorphic processes on $\Bil$ are
Markov-regular.
Now Markov-regularity of the process
corresponds to norm continuity of the convolution semigroup given by
\[
\lambda_t (F) =
\int_{\Omega} F \circ X_t\, d\, \Prob  \;\;
\]
($F \in \Bil, t \geq 0$).
Note that the usual notion of weak continuity for this semigroup corresponds,
in the algebraic formulation, to pointwise continuity of the Markov semigroup.
We therefore obtain the following result.

\begin{propn}
Let $X$ be a L\'evy processes on a compact semigroup with identity $H$. Suppose that as a topological
space $H$ is normal.
Then $X$ is equivalent
to a Markov-regular *-homomorphic QS convolution cocycle
on $\Bil$ if and only if it satisfies the following condition\tu{:}
\begin{equation}
\Prob \big(\{X_t = e \} \big) \to  1 \text{ as } t\to 0.
\label{nor}\end{equation}
\end{propn}
\begin{proof}
It is easily seen that condition (\ref{nor}) implies the existence of a bounded
generator $\gamma: \Bil \to \bc$ from which the process can be reconstructed. The other
direction can be seen by considering the Markov semigroup of a given QS convolution cocycle and
judiciously choosing continuous functions on $H$ with values in $[0,1]$, which are equal $1$ at
$e$ and vanish outside of some neighbourhood of the identity element $e$.
\end{proof}
Processes satisfying (\ref{nor}) were investigated for example in \cite{Grenander}. They are called
homogenous processes of discontinuous type and their laws are compound Poisson
distributions (\cite{Grenander}, Theorem 2.3.5).

\subsection{Cocommutative case: group algebras}

Let $\Gamma$ be a discrete group. Denote by $\Bil=C^*(\Gamma)$ the enveloping $C^*$-algebra
of the Banach algebra $l^1(\Gamma)$ (\cite{groupalgebras}), called the universal
(or full) $C^*$-algebra of
$\Gamma$. By construction (the algebra of functions on $\Gamma$
with finite support being dense in $\Bil$), there is a universal unitary representation
$L:\Gamma \to \Bil$ such that
$\Blg:= \text{Lin}\{L_g : g \in \Gamma\}$
is dense in $\Bil$.  Due to
universality the mappings $\Com$ and $\Cou$ defined on the image of $L$ by
\[
\Com(L_g) = L_g \ot L_g \text{ and } \Cou(L_g) =1,
\]
extend to *-homomorphisms on $\Bil$. It is easily checked that $\Bil$,
equipped with the resulting comultiplication and counit, becomes a cocommutative
$C^*$-bialgebra.

\begin{thm}
Let $\Bil = C^*(\Gamma)$ for a discrete group $\Gamma$. Then
\begin{equation}\label{W and l}
W(t,g) = l_t(L_g) \quad (g\in\Gamma, t\geq 0)
\end{equation}
defines a bijective correspondence between unital *-homomorphic
QS convolution cocycles on the $C^*$-bialgebra $\Bil$ and
maps $W:\br_+ \times\Gamma \to B(\FFock)$ satisfying the following conditions:
\begin{rlist}
\item
for each $g \in \Gamma$ the family $\{W(t,g): t \geq 0\}$ is a left
QS operator cocycle\tu{;}
\item
for each $t \geq 0$ the family $\{W(t,g) :g \in \Gamma\}$ is a unitary representation of
$\Gamma$ on $\FFock$.
\end{rlist}
\label{corp}
\end{thm}

\begin{proof}
Let $l \in \cProc(\Bil)$ is a *-homomorphic QS convolution
cocycle and define a map $W:\Rplus\times\Gamma\to B(\FFock)$
by~\eqref{W and l}.
Then, for all $g,h\in\Gamma$ and $s,t\geq 0$,
\begin{align*}
&l_{s+t} (L_g) = \big(l_s \ot (\sigma_t \circ l_s)\big)(\Com L_g)
= l_s(L_g) \ot \sigma_s(l_t (L_g)) = W(s,g) \ot \sigma_s (W(t,g)), \\
&l_t(L_g) l_t(L_h) = l_t (L_g L_h) = l_t (L_{gh}) = W(t,gh), \\
&l_t(L_g)^* = l_t(L_g^*) = l_t (L_{g^{-1}}) = W(t, g^{-1}),  \\
&l_t(L_e) = l_t (\ida) = I_{\FFock} \text{ and } \\
&l_0(L_g) = I_{\FFock},
\end{align*}
so $W$ satisfies {\rm (i) and (ii)}. Conversely, suppose that
$W:\br_+ \times\Gamma \to B(\FFock)$ is a map satisfying conditions
{\rm (i) and (ii)}.
Due to universality there are maps $l_t:\Bil\to B(\FFock)$, $t\geq 0$,
satisfying~\eqref{W and l}.
The properties of $W$ imply that they are unital *-homomorphisms
and that they satisfy
\[
l_0 (x) =
\Cou (x)I_{\FFock} \text{ and } l_{s+t} (x) =
\big(l_s \otimes (\sigma_s \circ l_t)\big)(\Com x)
\]
for $s,t\geq 0$ and $x \in \Blg$.
Continuity ensures that these remain valid for $x \in \Bil$ and so the
result follows.
\end{proof}

\noindent
On the level of stochastic generators the above correspondence takes the
following form.

\begin{thm} \label{corg}
Let $\Blg:= \Lin\{L_g : g \in \Gamma\}$ for a discrete group
$\Gamma$. Then
\[
\psi_g = \varphi (L_g), \quad g\in\Gamma,
\]
determines a bijective correspondence between
maps $\varphi \in L( \Blg; B(\kilhat))$ satisfying
\begin{equation}
\varphi (ab) = \varphi (a) \Cou(b) + \Cou(a) \varphi(b) + \varphi(a)
       \Projk \varphi(b), \;\; \varphi (a)^* = \varphi(a^*),
       \;\;\varphi(1) = 0, \label{pier}
\end{equation}
and maps $\psi : \Gamma \to B(\kilhat)$ satisfying
\begin{equation}    \label{drug}
 \psi_{gh} = \psi_g + \psi_h + \psi_g \Projk \psi_h,  \;\;
(\psi_g)^* = \psi_{g^{-1}},\;\; \psi_e = 0;
\end{equation}
\end{thm}

\begin{proof} Elementary calculation.
 \end{proof}

\begin{rems}
Identities (\ref{drug}) may be considered as a special (time-independent)
case of formulae (4.2-4) in \cite{flows}. They are equivalent to $\psi$
having the block matrix form
\begin{equation}\label{lambda xi U}
\psi_g =
\begin{bmatrix}
i \lambda_g - \frac{1}{2} \|\xi_g\|^2 & - \langle \xi_g |U_g \cr
 |\xi_g\rangle & U_g - I_{\kil} \cr
\end{bmatrix},
\end{equation}
for a unitary representation $U$ of $\Gamma$ on $\kilhat$
and maps $\lambda: \Gamma \to \Real$ and $\xi: \Gamma \to \kil$
satisfying
\[
\xi_{gh} = \xi_g + U_g \xi_h
\text{ and }
\lambda_{gh} =
\lambda_g + \lambda_h - {\rm Im} \langle \xi_g, U_g \xi_h \rangle.
\]
\end{rems}

Note that, according to Theorem 6.3 of \cite{LSqscc1}, each map
$\varphi\in L( \Blg; B(\kilhat))$ satisfying (\ref{pier})
generates a unital, real and weakly multiplicative QS convolution cocycle
$l^{\varphi}$ on $\Blg$. The process $l^{\varphi}$ continuously
extends to a *-homomorphic QS convolution cocycle on $\Bil$ (see Lemma
\ref{extend} below).
On the other hand, given a map $\psi: \Gamma\to B(\kilhat)$ satisfying
(\ref{drug}), for each fixed $g \in \Gamma$ the unique (weakly regular,
weak) solution of the operator QS differential equation
\[
X_0 = I_{\FFock}, \;\;\; dX_t = X_t  d\Lambda_L(t),
\]
where $L=\psi_g$, is a unitary left QS  cocycle $W^g$ (\cite{LWjfa}).
The map $W:\br_+ \times \Gamma$ given by
$W(t,g)= W_t^g$ satisfies the conditions of Theorem \ref{corp}. One can
easily see that
the correspondences described in Theorems \ref{corp} and \ref{corg} are consistent with this
construction.

\begin{propn}
A unital *-homomorphic QS convolution cocycle $l$ on $\Blg$
is equal to $l^{\varphi}$ for some $\varphi\in L( \Blg; B(\kilhat))$ if and only
if it is pointwise weakly measurable.
\end{propn}
\begin{proof}
One direction is trivial. For the other consider the unitary cocycles $\{W(\cdot, g):g \in \Gamma\}$ associated with $l$ by
Theorem \ref{corp}. Theorem 6.7 of \cite{LWjfa} implies that each of these cocycles is
stochastically generated (as it is weakly measurable). Denoting the respective generators by
$\psi_g$ one can see that the map $\psi:\Gamma \to B(\kil)$ so obtained
satisfies the conditions
(\ref{drug}). The desired conclusion therefore follows from Theorem \ref{corg}
and the subsequent discussion.
\end{proof}

If a *-homomorphic QS convolution cocycle $l$ on $\Bil$ is
Markov-regular, the automatic implementedness of its
stochastic generator $\varphi$ (Theorem~\ref{*homchar}) implies in particular that
the triple $(\lambda, \xi, U)$ corresponding
to $\varphi$ by~\eqref{lambda xi U} and~Theorem \ref{corg} must also be
\emph{implemented}, in the
following sense: there is a vector $\eta \in \kil$ such that
\[
\xi_g =
U_g \eta - \eta \text{ and } \lambda_g =
{\rm Im} \langle \eta, U_g \eta \rangle,
\quad g\in G.
\]
In the language of group cohomology, the first order cocycle $\xi$ is a
coboundary. In this connection, see~\cite{PaS}.

Elements of a $C^*$-bialgebra $\blg$ are called \emph{group-like}
when they satisfy $\Com b = b \ot b$, as the $L_g$'s do. On such elements
the solution $\big(k_t(b)\big)_{t\geq 0}$,
of the mapping QS differential equation~\eqref{OScoalgQSDE},
is given by the solution of the operator QS differential equation
\[
dX_t = X_t d \Lambda_L (t), \;\;\; X_0 = I_{\FFock},
\]
where $L= \varphi(b) \in B(\kilhat)$. For more on this we refer to Section
4.1 of \cite{Schurmann}.

\subsection{Full compact quantum groups}

A concept of compact quantum groups was introduced by Woronowicz, in \cite{wor1}. For our
purposes it is most convenient to adopt the following definition:

\begin{defn}[\cite{wor2}]
A \emph{compact quantum group} is a pair $(\Bil, \Com)$, where $\Bil$ is a unital $C^*$-algebra,
and $\Com:\Bil \to \Bil \ot \Bil$ is a unital, *-homomorphic map which is
coassociative and
satisfies the quantum cancellation properties:
\[ \overline{\Lin}((1\ot \Bil)\Com(\Bil) ) = \overline{\Lin}((\Bil \ot 1)\Com(\Bil) )
= \Bil \ot \Bil. \]
\end{defn}

For the concept of Hopf *-algebras and their unitary
corepresentations, as well as unitary corepresentations of compact quantum
groups, we refer the reader to~\cite{KlimykSchmudgen}.
For our purposes it is sufficient to note the facts contained in the
following theorem.

\begin{thm} [\cite{wor2}]
Let $\Bil$ be a compact quantum group and let $\Blg$ denote the linear
span of the matrix coefficients of irreducible unitary corepresentations
of $\Bil$. Then $\Blg$ is a dense *-subalgebra of $\Bil$,
the coproduct of $\Bil$ restricts to an algebraic coproduct
$\Com_0$ on $\Blg$ and
there is a natural counit $\Cou$ and coinverse $\Skorohod$ on $\Blg$ which
makes it a Hopf *-algebra.
\end{thm}

\begin{rem}[\cite{BMT}]
In the above theorem
$(\Blg,\Com_0, \Cou, \Skorohod)$ is the
unique dense Hopf *-subalgebra of $\Bil$, in the following sense: if
$(\Blg',\Com'_0, \Cou', \Skorohod')$
is a Hopf *-algebra, in which
$\Blg'$ is a dense *-subalgebra of $\Bil$ and
the coproduct of $\Bil$ restricts to the algebraic coproduct
$\Com'_0$ on $\Blg'$, then $(\Blg',\Com'_0, \Cou', \Skorohod')$ equals
$(\Blg,\Com_0, \Cou, \Skorohod)$.
\end{rem}

The Hopf  *-algebra arising here is
called the \emph{associated Hopf  *-algebra} of $(\Bil, \Com)$.
When $\Bil = C(G)$ for a compact group $G$, $\Blg$ is the algebra of all
matrix coefficients of unitary representations of $G$; when $\Bil$ is
the universal $C^*$-algebra of a discrete group $\Gamma$, $\Blg=\text{Lin}
\{L_g: g\in \Gamma\}$ (see the beginning of the previous subsection).
Dijkhuizen and Koornwinder observed that the Hopf *-algebras arising in
this way have intrinsic algebraic structure.

\begin{defn}
A Hopf  *-algebra $\Blg$ is called a
\emph{CQG algebra} if it is the linear span of all matrix
elements of its finite dimensional unitary corepresentations.
\end{defn}

\begin{thm}[\cite{DiK}]
Each Hopf  *-algebra associated with a compact
quantum group is a CQG algebra.
Conversely, if  $\Blg$ is a CQG algebra then
\begin{equation}
\label{norm on a CQG}
\|x\|:=
\sup \big\{\|\pi(x)\|:
\pi \text{ is a *-representation of } \Blg
\text{ on a Hilbert space}\big\}
\end{equation}
defines a $C^*$-norm on $\Blg$ and
the completion of $\Blg$ with respect to this norm is a compact
quantum group whose comultiplication extends that of
$\Blg$.
\end{thm}

\noindent
The compact quantum group obtained from a Hopf *-algebra $\Blg$
in this theorem is called its
\emph{universal compact quantum group} and is
denoted $\Blg_{\rm u}$.

For later use note the following extension of Lemma 11.31 in
\cite{KlimykSchmudgen}:

\begin{lemma}              \label{extend}
Let $E$ be a dense subspace of a Hilbert space $\Hil$ and let $\Blg$ be a
CQG algebra.
Suppose that $\pi:\Blg \to \Opdagger(E)$ is real, unital and weakly
multiplicative.
Then $\pi$ is bounded-operator-valued and admits a continuous extension to
a unital *-homomorphism from $\Blg_u$ to $B(\Hil)$.
\end{lemma}

\begin{proof}
Let $[x_{i,j}]_{i,j=1}^n$ be any finite dimensional unitary
corepresentation of $\Blg$. Then, since
$\sum_{k=1}^n  x_{k,j}^* x_{k,j} = 1_{\Blg}$
for $j \in \{1, \cdots,n\}$,
\begin{align*}
\| \pi (x_{i,j}) \xi \|^2
\leq \sum_{k=1}^n \| \pi (x_{k,j}) \xi \|^2
&=\sum_{k=1}^n \la \pi(x_{k,j}) \xi , \pi(x_{k,j}) \xi \ra \\
&=\left\la \xi , \pi \left( \sum_{k=1}^n  x_{k,j}^* x_{k,j} \right)\xi
\right \ra =
\|\xi \|^2
\end{align*}
for $i,j \in \{1, \cdots,n\}$ and $\xi \in E$.
This implies that, for each $x\in \Blg$,
$\pi(x)$ is bounded --- let $\pi_1 (a)$ denote its continuous extension to
a bounded operator on $\Hil$.
The resulting map $\pi_1 : \Blg \to B(\Hil)$ is then a unital
*-homomorphism, moreover it is clearly contractive
with respect to the canonical norm on $\Blg$, given
by~\eqref{norm on a CQG}; the result follows.
\end{proof}

\begin{defn}
A compact quantum group $(\Bil, \Com)$ is called \emph{full}
 if the $C^*$-norm
it induces on its associated CQG algebra $\Blg$
coincides with its canonical norm defined in \eqref{norm on a CQG}
--- equivalently, if $\Bil$ is *-isomorphic to $\Blg_{\rm u}$.
\end{defn}

The notion of full compact quantum groups was introduced in \cite{BMT} and in \cite{Baaj}
(in the first paper they were called universal compact quantum groups).  It is very relevant for our
 context, as the above facts imply the following

\begin{propn}
Let $\Bil$ be a full compact quantum group with associated Hopf *-algebra
$\Blg$. Then $\Bil$ is
a $C^*$-bialgebra whose counit is the continuous extension of the counit on $\Blg$.
Moreover restriction induces a bijective correspondence
between unital, *-homomorphic QS convolution cocycles on $\Bil$ and
unital, real, weakly multiplicative QS convolution cocycles \tu{(}in the
sense of~\cite{LSqscc1}\tu{)} on $\Blg$.
\end{propn}

Both families of examples described in the previous two subsections, namely algebras of continuous
functions on compact groups and full $C^*$-algebras of discrete groups, are full compact
quantum groups. Moreover most of the genuinely quantum (i.e. neither commutative nor
cocommutative) compact quantum groups considered in the literature also fall into this category,
including the queen of examples,  $SU_q(2)$.

Reconnecting further with our previous work, we obtain the following result.

\begin{thm}
Let $k\in\cProccb(\Bil)$ where $\Bil$ is a full compact quantum
group with associated Hopf *-algebra $\Blg$.
Then the following are equivalent\tu{:}
\begin{rlist}
\item
$k$ and $k^{\dagger}$ are pointwise H\"older-continuous QS convolution cocycles;
\item
$k|_{\Blg}=l^{\varphi}$ for some map $\varphi \in L(\Blg; B(\kilhat))$.
\end{rlist}
\end{thm}

\begin{proof}
One direction follows from the fact that $\Blg$ is an (algebraic) coalgebra and Theorem 5.8 of
\cite{LSqscc1}. The other is trivial.
\end{proof}

Specialising to *-homomorphic cocycles yields the following much stronger
result.

\begin{thm}
Let $k\in\cProc(\Bil :\ExpsD)$
where $\Bil$ is a full compact quantum group with asociated Hopf *-algebra $\Blg$ and $D$ is a dense
subspace of $\kil$.
Then the following are equivalent\tu{:}
\begin{rlist}
\item
$k$ is pointwise H\"older-continuous, unital and *-homomorphic (thus bounded) and $a \longmapsto k_t(a)$ defines a QS convolution cocycle\tu{;}
\item
$k$ is bounded and $k|_{\Blg}=l^{\varphi}$ for some
 $\varphi \in L(\Blg; \Opdagger(\Dhat))$ satisfying the structure
relations~\eqref{hom structure}.
\end{rlist}
\end{thm}

\begin{proof}  The implication (i)$\Rightarrow$(ii)
follows from the previous theorem and implication (i)$\Rightarrow$(ii) of
Theorem 6.3 of \cite{LSqscc1}
(note that it even yields $\varphi \in L(\Blg; \Opdagger(\kilhat)) =
L(\Blg; B(\kilhat))$.

Suppose conversely that (ii) holds. Theorem 6.3 of \cite{LSqscc1}
guarantees that $l=k|_{\Blg}$ is real, unital, and weakly multiplicative.
Lemma \ref{extend} shows that $l$ admits a continuous extension to a
*-homomorphic unital process $\Bil \to \bc$ defined on $\Exps_D$, which
must coincide with $k$. Application of the previous theorem therefore
completes the proof.
\end{proof}

The above theorem may be equivalently formulated in the following way.

\begin{thm}
Let $k\in\cProc\big(\Blg :\ExpsD\big)$ where
$\Blg$ is the Hopf *-algebra associated with  a full compact quantum group
$\Bil$ and $D$ is a dense
subspace of $\kil$.
Then the following are equivalent\tu{:}
\begin{rlist}
\item
$k$ extends to a  pointwise H\"older-continuous, unital, *-homomorphic
QS convolution cocycle on $\Bil$\tu{;}
\item
$k=l^{\varphi}$ for some
 $\varphi \in L\big(\Blg; \Opdagger(\Dhat)\big)$ satisfying the structure
relations~\eqref{hom structure}.
\end{rlist}
\end{thm}

\begin{rem}
In the course of the proof of the previous theorem it was established that each map $\varphi$ defined on
a CQG algebra $\Blg$ with values in $\Opdagger(\Dhat)$ satisfying the
conditions
(\ref{hom structure}) must be bounded-operator-valued. We stress however, that
$\varphi$ need not extend continuously to $\Bil$
(for examples see \cite{SchSk}). On the other hand if $\varphi$ is continuous,
then it is necessarily completely bounded.
\end{rem}

\renewcommand{\theequation}{\Alph{section}.\arabic{equation}}

\renewcommand{\thesection}{\Alph{section}}

\setcounter{section}{0}

\renewcommand{\thepropn}{\Alph{section}.\arabic{propn}}

\section*{Appendix: $(\pi',\pi)$-derivations and $\chi$-structure maps}

\setcounter{section}{1}
\setcounter{propn}{0}

In this appendix we give an extension of the
innerness theorem of Christensen, for completely bounded derivations on a
$C^*$-algebra, to $(\pi',\pi)$-derivations, and prove automatic complete
boundedness for $(\pi,\chi)$-derivations, when $\chi$ is a character.
These are then applied to prove the innerness of what we call $\chi$-structure maps.
We first recall the relevant theorems on derivations.

\begin{thm} [\cite{Sakai}, \cite{Ringrose}] \label{Ringrose}
Let $\delta :\alg\to\Xil$ be a derivation from a $C^*$-algebra $\alg$
into a Banach $\alg$-bimodule. Then $\delta$ is bounded.
\end{thm}

\begin{thm} [\cite{Der}]    \label{Christens}
Let $\alg$ be a $C^*$-algebra in $B(\hil)$ and let $\delta: \alg \to
B(\hil)$ be a derivation. If $\delta$ is completely bounded then it is
inner\tu{:} there is $R\in B(\hil)$ such that $\delta(a) = aR-Ra,  a \in \alg$.
\end{thm}

A simple proof of the first theorem in the case $\Xil= \alg$ (Sakai's
Theorem), due to Kishimoto, may be found in~\cite{SakaiBook}, and a good reference
for the second, along with connections to not-necessarily-involutive
homomorphisms between $C^*$-algebras, is~\cite{Pisiersim}.
We are interested in the particular class of Banach $\Al$ bimodule-valued
derivations captured by the following definition.

\begin{defn}
Let $\alg$ be a $C^*$-algebra with representations
$(\pi, \hil)$ and $(\pi', \hil')$.
A map $\delta: \alg \to B(\hil;\hil')$ is called a
$(\pi',\pi)$-\emph{derivation}
if it satisfies
\[
\delta(ab)  = \delta(a) \pi(b) +
\pi'(a) \delta(b);
\]
it is \emph{inner} if it is implemented by an operator
$T \in B(\hil;\hil')$ in the sense that
\[
\delta : a \to \pi'(a) T - T \pi(a).
\]
\end{defn}

\begin{thm} \label{RingChristen}
Let $\alg$ be a $C^*$-algebra with representations $(\pi, \hil)$ and
$(\pi', \hil')$, and let $\delta: \alg \to B(\hil;\hil')$ be a
completely bounded $(\pi',\pi)$-derivation.
Then $\delta$ is inner.
\end{thm}

\begin{proof}
Let $(\rho, \Kil)$ be a faithful representation of $\alg$ and set
$\Hil = \hil \oplus \hil' \oplus \Kil$ and
$\wt{\alg} = \wt{\pi}(\alg)$ where $\wt{\pi}$ is the faithful
representation
$\pi \oplus \pi' \oplus \rho$. Then $\wt{\alg}$ is a  $C^*$-subalgebra
of $B(\Hil)$ and
it is easily verified that
\[
\wt{\pi}(a)
\mapsto
\begin{bmatrix}
0 & & \cr
\delta(a) & 0 & \cr
& & 0
\end{bmatrix}
\]
defines a derivation
$\wt{\delta}:\wt{\alg} \to B(\Hil)$. It is also clear that $\wt{\delta}$
is completely bounded if and only if $\delta$ is.
Moreover, if $\wt{\delta}$ is inner then the
$(\pi',\pi)$-derivation $\delta$ is implemented by
$S_{21}\in B(\hil;\hil')$ for any operator $S=[S_{ij}]\in B(\Hil)$
implementing the derivation $\wt{\delta}$. The result therefore follows
from Theorem~\ref{Christens}.
\end{proof}

\begin{thm}
\label{cbderiv}
Let $\alg$ be a $C^*$-algebra with representation $(\pi, \hil)$ and
character $\chi$, and let $\delta: \alg \to |\hil\ra$ be a
$(\pi, \chi)$-derivation. Then $\delta$ is inner.
\end{thm}

\begin{proof}
Without loss of generality we may suppose that the $C^*$-algebra $\Al$ and
representation $\pi$ are both unital; if necessary by extending $\pi$,
$\chi$ and $\delta$ to the unitisation of $\Al$ in the following natural
way:
\[
(a,z)\mapsto \pi (a) + zI_{\hil}, \ \
(a,z)\mapsto \chi (a) + z \ \text{ and } \
(a,z)\mapsto \delta (a).
\]
By Theorem~\ref{Ringrose}, $\delta$ is bounded.
Let $\Ao = \Ker \chi$ and let $\psi : \Al \rightarrow \Ao$ be the
projection
$a \mapsto a - \chi (a) 1$. Then $\Ao$ is a $C^*$-subalgebra of $\Al$,
$\psi$ is
completely bounded and $\delta = \wt{\delta} \circ \psi$, where
$\wt{\delta} = \delta |_{\Ao}$.
Therefore, by the previous theorem, it suffices to show that $\wt{\delta}$
is completely bounded.
Now $\wt{\delta}(ab) = \pi (a) \delta (b)$ for all $a,b \in \Ao$.  Since
$\delta$ is bounded this implies that
\[
\wt{\delta}^{(n)} (A) =
\lim_{\lambda} \pi^{(n)} (A) \big( \delta (e_\lambda ) \ot I_n \big)
\]
($n\in\Nat, A \in \Mat_n (\Ao )$),
for any $C^*$-approximate identity $(e_\lambda)$ for $\Ao$, and so
 $\big\| \wt{\delta}^{(n)} \big\| \leq  \| \delta \|$. The result follows.
\end{proof}

We note two consequences; the first is used in~\cite{Sthesis}.

\begin{cor} \label{disap}
Let $\alg$ be a $C^*$-algebra with characters (i.e.\ nonzero multiplicative linear functionals) $\chi$ and $\chi'$.
Then every $(\chi',\chi)$ derivation on $\alg$ vanishes.
\end{cor}

For the second the following definitions are convenient. If $\Al$ is a
$C^*$-algebra with character $\chi$, then a
$\chi$-\emph{structure map on} $\Al$
is a linear map $\varphi : \Al \to B(\Comp \op \hil)$, for some Hilbert
space $\hil$, satisfying
\begin{equation} \label{chi structure relation}
\varphi (a^*b) = \varphi(a)^* \chi (b) + \chi(a)^* \varphi(b) + \varphi (a)^* \Delta \varphi (b)
\end{equation}
where $\Delta := \left[ \begin{smallmatrix} 0 & \\ & I_{\hil} \end{smallmatrix} \right]$. For any
$C^*$-representation $(\pi , \hil)$ and vector $\xi \in \hil$,
\[
a \mapsto \begin{bmatrix} \la \xi | \\ I_{\hil} \end{bmatrix}
\big( \pi (a) - \chi (a) I_{\hil} \big) \begin{bmatrix} | \xi \ra & I_{\hil} \end{bmatrix}
\]
defines a $\chi$-structure map. Such $\chi$-structure maps are said to be \emph{implemented}.
Thus implementation involves a pair $(\pi, \xi)$. Note
that implemented $\chi$-structure maps are completely bounded.

\begin{thm} \label{A6}
Let $\Al$ be a $C^*$-algebra with character $\chi$ and let
$\varphi$ be a $\chi$-structure map on $\Al$. Then $\varphi$ is implemented.
\end{thm}

\begin{proof}
Without loss of generality we may suppose that $\Al$ is unital, since
otherwise (invoking the reality of $\varphi$) the prescriptions
\[
(a,z)\mapsto\chi (a) + z, \ \text{ respectively } \
(a,z)\mapsto\varphi (a),
\]
extend $\chi$ and $\varphi$ to the unitisation of $\Al$, maintaining
the $\chi$-structure relation~\eqref{chi structure relation}.
Now the $\chi$-structure relation is equivalent to $\varphi$ having
block matrix form
\begin{equation} \label{the new 6.6}
\begin{bmatrix} \lambda & \delta^\dagger \\ \delta & \nu \end{bmatrix}
\end{equation}
where $\nu = \pi - \iota_{\kil} \circ \chi$ for a *-homomorphism
$\pi :\Al \to B(\hil)$, $\delta$
is a $(\pi , \chi )$-derivation and the linear functional $\lambda$ satisfies
\[
\lambda (a^*b) = \lambda(a)^* \chi(b) + \chi (a)^* \lambda (b) + \delta(a)^* \delta (b)
\]
$(a,b \in \Al)$ --- in particular, $\lambda$ is real and satisfies
\begin{equation} \label{lambda 1 A0}
\lambda (1) = - \delta (1)^* \delta(1) \text{ and }
\lambda (a^*b) = \delta(a)^* \delta(b) \text{ for } a,b \in \Al_0 ,
\end{equation}
where $\Al_0 = \Ker \chi$. By Theorem~\ref{cbderiv}, there is a vector
$\xi \in \kil$ such that $\delta (a) = \nu (a) | \xi \ra$. Now define a
bounded linear functional $\wt{\lambda}$ on $\Al$
by $\wt{\lambda} (a) = \big\la \xi , \nu (a) \xi \big\ra$.
It is easily checked that
$\wt{\lambda}$ also satisfies~\eqref{lambda 1 A0}, thus
$\wt{\lambda}$ agrees with $\lambda$ on $\Al_{00} + \Comp 1_{\Al}$ where
$\Al_{00} = \Lin \{ a^*b :a,b \in \Al_0 \}$. But $\Al_{00}$ is dense in
$\Al_0$ and $\Al=\Al_0 \op \Comp 1_{\Al}$,
so $\wt{\lambda}$ equals $\lambda$. The result follows.
\end{proof}

\end{document}